\documentclass{cmslatex}

\usepackage{ulem}
\usepackage{epsfig}
\usepackage{graphics} 
\DeclareGraphicsExtensions{.pdf}
\graphicspath{
{/u/cermics/r/boyaval/Desktop/stochRB/reducvar_controlvar/CMS/fig_pdf/}
}\usepackage{amssymb}
\usepackage{amsmath}
\usepackage[latin1]{inputenc}
\usepackage[T1]{fontenc}

\newlength{\mywidth}
\renewcommand{\theequation}{\arabic{section}.\arabic{equation}} 
\newtheorem{theo}{Theorem}[section]
\newtheorem{prop}[theo]{Proposition}

\newtheorem{rem}[theo]{Remark}
\def\be{\begin{equation} \displaystyle}
\def\ee{\end{equation} }

{ \begin{list}%
	{-}
	{\setlength{\labelwidth}{5pt}%
	 \setlength{\leftmargin}{5pt}%
	 \setlength{\itemsep}{\parsep}
	}
}
{ \end{list} }

\newcommand{\uu}[1]  {{\boldsymbol #1} }
\newcommand{\argmax}{\operatornamewithlimits{arg\,max}}

\def\P{\mathbb{P}} 
\newcommand{\E}[1]{\mathbb{E}\left(#1\right)}
\newcommand{\Var}[1]{{\rm\mathbb{V}ar}\left(#1\right)}
\newcommand{\Cov}[1]{{\rm\mathbb{C}ov}\left(#1\right)}
\def\tr{\mathrm{tr}} 
\newcommand{\Span}[1]{\mathbf{Span}\left(#1\right)}
\def\I{\mbox{Id}} 
\def\F{\mathcal{F}} 
\def\N{\mathbb{N}}
\def\R{\mathbb{R}}

\def\dt{\Delta t}

\newcommand{\comment}[1]{ }

\renewcommand{\emph}[1]{\textbf{#1}}

\begin{document}

\title{
A Variance Reduction Method\\
for Parametrized Stochastic Differential Equations\\
using the Reduced Basis Paradigm
}

\author{
S\'ebastien Boyaval
\thanks{
CERMICS, Ecole des Ponts ParisTech (Universit\'e Paris-Est), Cit\'e Descartes,
77455 Marne-la-Vall\'ee Cedex 2, France, \ and \ 
MICMAC team-project, INRIA, Domaine de Voluceau, BP. 105  Rocquencourt
78153 Le Chesnay Cedex, France.
}
\and 
Tony Leli\`evre$\ {}^{\dagger}$
}



\maketitle

\renewcommand{\thefootnote}{\*}
\footnotetext{Corresponding author: S.~Boyaval. \\
Email address: boyaval@cermics.enpc.fr \\
URL home page: http://cermics.enpc.fr/$\sim$boyaval/home.html \\
Tel: + 33 1 64 15 35 79 - Fax: + 33 1 64 15 35 86
}

\renewcommand{\thefootnote}{\arabic{footnote}}
\setcounter{footnote}{0}



\begin{abstract}
In this work, we develop a reduced-basis approach for the efficient computation of 
parametrized expected values,
for \textit{a large number} of parameter values, 
using the control variate method to reduce the variance.
Two algorithms are proposed to compute \textit{online}, 
through a cheap reduced-basis approximation,
the control variates 
for the computation of a large number of expectations of a functional of a parametrized It\^o stochastic process
(solution to a parametrized stochastic differential equation).
For each algorithm, a reduced basis of control variates is pre-computed \textit{offline},
following a so-called \textit{greedy} procedure, 
which minimizes the variance among a trial sample of the output parametrized expectations.
Numerical results in situations relevant to practical applications
(calibration of volatility in option pricing, and
parameter-driven evolution of a vector field following a Langevin equation from kinetic theory)
illustrate the efficiency of the method.
\end{abstract}

\section{Introduction}

This article develops a general variance reduction method for the \emph{many-query} context 
where a large number of \emph{Monte-Carlo} estimations of the expectation $\E{Z^{\lambda}}$ of a functional 
\begin{equation}\label{functional}
Z^{\lambda}=g^\lambda(X_T^{\lambda})-\int_0^T f^\lambda(s,X_s^{\lambda})\, ds
\end{equation}
of the solutions $\left( X_t^{\lambda} , t\in[0,T] \right)$ 
to the \emph{stochastic} \emph{differential} \emph{equations} (SDEs):
\begin{equation}\label{eq:pb}
X_t^{\lambda} = x + \int_0^t b^\lambda(s,X_s^{\lambda})\, ds 
	+ \int_0^t  \sigma^{\lambda}(s,X_s^{\lambda}) dB_s
\end{equation}
parametrized by $\lambda \in \Lambda$ have to be computed for many values of the parameter $\lambda$.

Such many-query contexts are encountered in finance for instance, 
where pricing options often necessitates to compute the price $\E{Z^{\lambda}}$ of an option
with spot price $X_t^{\lambda}$ at time $t$ in order to \textit{calibrate} the local volatility
$\sigma^{\lambda}$ as a function of a (multi-dimensional) parameter $\lambda$
(that is minimize over $\lambda$, after many iterations of some optimization algorithm,
the difference between observed statistical data with the model prediction).
Another context for application is molecular simulation,
for instance micro-macro models in rheology, 
where the mechanical properties of a flowing viscoelastic fluid are determined 
from the coupled evolution of a non-Newtonian stress tensor field $\E{Z^{\lambda}}$
due to the presence of many polymers with configuration $X_t^{\lambda}$ 
in the fluid with instantaneous velocity gradient field $\lambda$.
Typically, segregated numerical schemes are used: compute $X_t^{\lambda}$ for a fixed field $\lambda$, 
and then compute $\lambda$ for a fixed field $\E{Z^{\lambda}}$.
Such tasks are known to be computationally demanding
and the use of different \emph{variance} \emph{reduction} techniques to alleviate the cost of Monte-Carlo computations 
in those fields is very common
(see~\cite{arouna-03,melchior-ottinger-95,ottinger-vandenbrule-hulsen-97,bonvin-picasso-99} for instance).

In the following, we focus on one particular variance reduction strategy 
termed the \emph{control} \emph{variate}
method~\cite{hammersley-handscomb-64,newton-94,milstein-tretyakov-06}.
More precisely, we propose new approaches in the context of the computation of $\E{Z^{\lambda}}$
for a large number of parameter values $\lambda$, with the control variate method.
In these approaches,
the control variates are computed through a \emph{reduced-basis} method
whose principle is related to 
the reduced-basis method~\cite{machiels-maday-patera-01,
maday-patera-turinici-02,patera-rozza-07,boyaval-08,boyaval-lebris-maday-nguyen-patera-08}
previously developed to efficiently solve parametrized partial differential equations (PDEs).
Following the reduced-basis paradigm,
a small-dimensional vector basis is first built \textit{offline}
to span a good linear approximation space for
a large trial sample of the $\lambda$-parametrized control variates,
and then used \textit{online} to compute control variates at any parameter value.
The offline computations are typically expensive, but done once for all. 
Consequently, it is expected that the online computations 
(namely, approximations of $\E{Z^{\lambda}}$ for many values of $\lambda$) are very cheap, 
using the small-dimensional vector basis built offline for \textit{efficiently} computing control variates online.
Of course, such reduced-basis approaches can only be \textit{efficient} insofar as:
\begin{enumerate}
 \item online computations (of one output $\E{Z^{\lambda}}$ for one parameter value $\lambda$)
are significantly cheaper using the reduced-basis approach than without, and
 \item the amount of outputs $\E{Z^{\lambda}}$ to be computed online 
(for many different parameter values $\lambda$) is sufficient to compensate for 
the (expensive) offline computations (needed to build the reduced basis).
\end{enumerate}
In this work, we will study numerically 
how the variance is reduced in two examples
using control variates built with two different approaches. 
\comment{
As a consequence,
it is possible to roughly estimate regimes of efficiency 
for our two different reduced-basis approaches,
that is when usual methods for the evaluation of the same outputs 
would need more
computations to reach the same online variance as our method,
taking into account the offline effort.
}

The usual reduced-basis approach for parametrized PDEs
also traditionally focuses on the certification of the reduction (in the parametrized solution manifold)
by estimating \textit{a posteriori} the error between approximations obtained before/after reduction
for some \textit{output} which is a functional of the PDE solution.
Our reduced-basis approach for the parametrized control variate method
can also be cast into a goal-oriented framework similar to the traditional reduced basis method.
One can take the expectation $\E{Z^{\lambda}}$ as the reduced-basis output,
while the empirically estimated variance ${\rm Var_M}\left(Z^{\lambda}\right)$
serves as a computable (statistical) error indicator for the Monte-Carlo
approximations 
${\rm E_M}\left(Z^{\lambda}\right)$ 
of $\E{Z^{\lambda}}$ \textit{in the limit of large $M$} through the Central Limit Theorem 
(see error bound~\eqref{eq:clt} in Section~\ref{sec:math}).

In the next Section~\ref{sec:var_red}, the variance reduction issue and the control variate method are introduced,
as well as the principles of our reduced-basis approaches for the computation of parametrized control variates.
The Section~\ref{sec:approximate} exposes details about the algorithms
which are numerically applied to 
test problems in the last Section~\ref{sec:numerics}.

The numerical simulations show good performance of the method for 
the two test problems corresponding
to the applications mentionned above:
a scalar SDE with (multi-dimensional) parametrized diffusion
(corresponding to the calibration of a local volatility in option pricing),
and a vector SDE with (multi-dimensional) parametrized drift
(for the parameter-driven evolution of a vector field following a Langevin equation from kinetic theory).
Using the control variate method with a $20$-dimensional reduced basis of (precomputed) control variates, 
the variance is approximatively divided
by a factor of $10^4$ in the mean for large test samples of parameter
in the applications we experiment here.
As a consequence,
our reduced-basis approaches allows 
to approximately divide the online computation time 
by a factor of $10^2$,
while maintaining the confidence intervals for the output expectation 
at the same value than without reduced basis.

This work intends to present a new numerical method and 
to demonstrate its interest on some relevant test cases.
We do not have, for the moment, a theoretical understanding of the method.
This is the subject of future works.

\section{The variance reduction issue and the control variate method}\label{sec:var_red}

\subsection{Mathematical preliminaries and the variance reduction issue}
\label{sec:math}

Let $\left( B_t \in \R^d , t\in[0,T] \right)$ be a $d$-dimensional standard Brownian motion 
(where $d$ is a positive integer)
on a complete probability space $(\Omega,\F,\P)$, endowed with a filtration $\left(\F_t,t\in[0,T]\right)$.
For any square-integrable random variables $X,Y$ on that probability space $(\Omega,\F,\P)$,
we respectively denote by $\E{X}$ and $\Var{X}$ the expected value and the variance of $X$ 
with respect to the probability measure $\P$,
and by $\Cov{X;Y}$ the covariance between $X$ and $Y$.

For every $\lambda \in \Lambda$ ($\Lambda$ being the set of parameter values),
the It\^o processes $\left(X_t^{\lambda}\in \R^d,t\in[0,T]\right)$ 
with deterministic initial condition $x\in \R^d$
are well defined as the solutions to the SDEs~\eqref{eq:pb} 
under suitable assumptions on $b^\lambda$ and $\sigma^\lambda$,
for instance provided $b^\lambda$ and $\sigma^\lambda$ satisfy Lipschitz and growth conditions~\cite{kloeden-platen-00}.
Let $(X_t^{\lambda})$ be solutions to the SDEs,
and $f^\lambda$, $g^\lambda$ be measurable functions such that 
$Z^{\lambda}$ is a well-defined integrable random variable ($Z^{\lambda}\in L^1_\P(\Omega)$).
Then, Kolmogorov's strong law of large numbers holds and, 
denoting by $Z_m^{\lambda}$ ($m=1,\ldots,M$) $M$ independent copies of the random variables $Z^{\lambda}$ 
(for all positive integer $M$),
the output expectation $\E{Z^{\lambda}}=\int_\Omega Z^{\lambda} d\P$ can be approximated (almost surely)
by Monte-Carlo estimations of the form:
\begin{equation} \label{eq:empirical_estimation}
 {\rm E}_M\left(Z^{\lambda}\right) := \frac1M \sum_{m=1}^M Z_m^{\lambda} 
\xrightarrow[M\to\infty]{\P-a.s.} \E{Z^{\lambda}} .
\end{equation}
Furthermore, assume that the random variable $Z^{\lambda}$ is square integrable 
($Z^{\lambda}\in L^2_\P(\Omega)$) with variance $\Var{Z^{\lambda}}$, 
then an asymptotic error bound for the convergence occuring in~\eqref{eq:empirical_estimation}
is given in probabilistic terms by the Central Limit Theorem as
confidence intervals: for all $a > 0$,
\begin{equation}
\label{eq:clt2}
\P\left( 
\left| {\rm E}_M\left(Z^{\lambda}\right) - \E{Z^{\lambda}} \right|
\le a\sqrt{\frac{\Var{Z^{\lambda}}}{M}} \right)
\underset{M\to\infty}{\longrightarrow} \int_{-a}^a \frac{e^{-x^2/2}}{\sqrt{2\pi}} dx  \,.
\end{equation}

In terms of the error bound~\eqref{eq:clt2},
an approximation ${\rm E}_M\left(Z^{\lambda}\right)$ of the output $\E{Z^{\lambda}}$
is thus all the better, for a given $M$, as the variance $\Var{Z^{\lambda}}$ is small.
In a many-query framework, 
the computation of approximations~\eqref{eq:empirical_estimation} for many outputs $\E{Z^{\lambda}}$ 
(corresponding to many queried values of the parameter $\lambda\in\Lambda$)
would then be all the faster as the variance $\Var{Z^{\lambda}}$ for some $\lambda\in\Lambda$
could be decreased from some knowledge acquired from the $\lambda\in\Lambda$ computed beforehand.
This typically defines a many-query setting with parametrized output
suitable for a reduced-basis approach similar to 
the reduced-basis method developped in a deterministic setting for parametrized PDEs.

In addition, the convergence~\eqref{eq:empirical_estimation} controlled
by the confidence intervals~\eqref{eq:clt2}
can be easily observed using computable {\it a posteriori} estimators.
Indeed, remember that since the random variable $Z^{\lambda}$ has a finite second moment,
then the strong law of large numbers also implies the following convergence:
\be \label{eq:variance-estimation}
{\rm Var}_M\left(Z^{\lambda}\right):=
{\rm E}_{M}\left( \left(Z^{\lambda}-{\rm E}_M(Z^{\lambda})\right)^2 \right) \xrightarrow[M\to\infty]{\P-a.s.} \Var{Z^{\lambda}} \ .
\ee
Combining the Central Limit Theorem  with Slutsky theorem
for the couple of Monte-Carlo estimators
$\left({\rm E}_M\left(Z^{\lambda}\right),{\rm Var}_M\left(Z^{\lambda}\right)\right)$
(see for instance~\cite{grimmett-stirzaker-01}, exercise 7.2.(26)),
we obtain a fully computable probabilistic (asymptotic) error bound 
for the Monte-Carlo approximation~\eqref{eq:empirical_estimation} 
of the output expectation: for all $a > 0$,
\be
\label{eq:clt}
\P\left( \left| \E{Z^{\lambda}} - {\rm E}_M\left(Z^{\lambda}\right) \right|
\le a\sqrt{\frac{{\rm Var}_M\left(Z^{\lambda}\right)}{M}} \right)
\overset{M\to\infty}{\longrightarrow} \int_{-a}^a \frac{e^{-x^2/2}}{\sqrt{2\pi}} dx  \,.  
\ee

It is exactly the purpose of variance reduction techniques 
to reduce the so-called \emph{statistical} error
appearing in the Monte-Carlo estimation of the output expectation $\E{Z^{\lambda}}$
through the error bound~\eqref{eq:clt2}.
And this is usually achieved in practice 
by using the ({\it a posteriori}) estimation~\eqref{eq:clt}.

\begin{rem}[SDE discretization and bias error in the output expectation]\label{rem:sde_discretization}
In practice, there is of course another source of error, 
coming from the time-discretizations of the SDE \eqref{eq:pb} 
and of the integral involved in the expression for~$Z^{\lambda}$.

In the following (for the numerical applications), 
we use the Euler-Maruyama numerical scheme
with discretizations $ 0=t_0<t_1<\dots<t_N=T $ ($N\in\N$) of the time interval $[0,T]$
to approximate the It\^{o} process $(X_t^{\lambda})$:
\begin{equation*}\label{eq:euler-maruyama}
\left\{\begin{array}{l}
\overline{X}{}_n^{\lambda} = \overline{X}{}_{n-1}^{\lambda} 
	+ |t_n-t_{n-1}|\ b^\lambda(t_{n-1},\overline{X}{}_{n-1}^{\lambda})  
	+ \sqrt{|t_n-t_{n-1}|}\sigma^{\lambda}(t_{n-1},\overline{X}{}_{n-1}^{\lambda})G_{n-1},
\\
\overline{X}{}_0^{\lambda} = x,
\end{array}\right.
\end{equation*}
where $\{G_{n},\ n=0,\ldots,N-1\}$ is a collection of $N$ independent $d$-dimensional 
normal centered Gaussian vectors.
It is well-known that such a scheme if of weak order one, so that
we have a bound for the \emph{bias} due to the approximation of the output expectation
$\E{Z^{\lambda}}$ by $\E{\overline{Z}{}^{\lambda}}$ (where $\overline{Z}{}^{\lambda}$
is a time-discrete approximation for $Z^{\lambda}$ computed 
from $(\overline{X}{}_n^{\lambda})$ 
with an appropriate discretization of 
the integral $\int_0^T f^{\lambda}(s,X_{s}^{\lambda}) ds$):
\begin{equation*}
\left|\E{\overline{Z}{}^{\lambda}}-\E{Z^{\lambda}}\right|
\underset{N\to\infty}{=} O\left(\underset{1\le n\le N}{\max}(|t_n-t_{n-1}|)\right) \,. \displaystyle
\end{equation*}
The approximation of the output $\E{Z^{\lambda}}$
by ${\rm E_M}(\overline{Z}{}^{\lambda})$ thus contains two types of errors:
\begin{itemize}
 \item first, a bias $\E{Z^{\lambda}-\overline{Z}{}^{\lambda}}$
 due to discretization errors in the numerical integration of the SDE~\eqref{eq:pb}
 and of the integral involved in $Z^{\lambda}$,
 \item second, a statistical error of order $\sqrt{{\Var{\overline{Z}{}^{\lambda}}}/{M}}$ 
 in the empirical Monte-Carlo estimation ${\rm E_M}(\overline{Z}{}^{\lambda})$ 
 of the expectation $\E{\overline{Z}{}^{\lambda}}$.
\end{itemize}
We focus here on the statistical error.
\end{rem}

\subsection{Variance reduction with the control variate method}
\label{sec:optimal}

The idea of control variate methods for the Monte-Carlo evaluation of $\E{Z^{\lambda}}$
is to find a so-called {\it control variate} $Y^{\lambda}$ (with $Y^{\lambda}\in L^2_\P(\Omega)$), 
and then to write:
$$ \E{Z^{\lambda}} = \E{Z^{\lambda}-Y^{\lambda}} + \E{Y^{\lambda}} \,,$$
where $\E{Y^{\lambda}}$ can be easily evaluated,
while the expectation $\E{Z^{\lambda}-Y^{\lambda}}$ is approximated by Monte-Carlo estimations 
that have a smaller statistical error than direct Monte-Carlo estimations of $\E{Z^{\lambda}}$. 
In the following, we will consider control variates $Y^{\lambda}$ such that 
$ \E{Z^{\lambda}} = \E{Z^{\lambda}-Y^{\lambda}} $, equivalently
$$\E{Y^{\lambda}}=0\,.$$ 
The control variate method will indeed be interesting if 
the statistical error of the Monte-Carlo estimations ${\rm E}_M(Z^{\lambda}-Y^{\lambda})$
is significantly smaller than the statistical error of the Monte-Carlo estimations ${\rm E}_M(Z^{\lambda})$. 
That is, considering the following error bound given by the Central
Limit Theorem: for all $a > 0$,
\begin{equation}
\label{eq:clt3}
\P\left( 
\left| {\rm E}_M\left(Z^{\lambda}-Y^{\lambda}\right)-\E{Z^{\lambda}} \right|
\le a\sqrt{\frac{\Var{Z^{\lambda}-Y^{\lambda}}}{M}}
\right)
\overset{M\to\infty}{\longrightarrow} \int_{-a}^a \frac{e^{-x^2/2}}{\sqrt{2\pi}} dx  \,,
\end{equation}
the Monte-Carlo estimations ${\rm E}_M(Z^{\lambda}-Y^{\lambda})$
will indeed be more accurate approximations of the expectations $\E{Z^{\lambda}}$
than the Monte-Carlo estimations ${\rm E}_M(Z^{\lambda})$ provided:
$$ \Var{Z^{\lambda}}\ge\Var{Z^{\lambda}-Y^{\lambda}} \,. $$

Clearly, the best possible control variate (in the sense of minimal variance) 
for a fixed parameter $\lambda\in\Lambda$ is:
\be\label{eq:best}
Y^{\lambda}=Z^{\lambda}-\E{Z^{\lambda}}\,,
\ee
since we then have $\Var{Z^{\lambda}-Y^{\lambda}}=0$.
Unfortunately, the result $\E{Z^{\lambda}}$ itself is necessary to compute $Y^{\lambda}$ as $Z^{\lambda}-\E{Z^{\lambda}}$.

In the following, we will need another representation of the best possible control variate $Z^{\lambda}-\E{Z^{\lambda}}$.
Under suitable assumptions on the coefficients $b^\lambda$ and $\sigma^\lambda$
(for well-posedness of the SDE), plus continuity and polynomial growth 
conditions on $f^{\lambda}$ and $g^{\lambda}$,
let us define $u^\lambda(t,y)$, for $(t,y)\in [0,T] \times \R^d$, as 
the unique solution $u^\lambda(t,y)\in C^1\left([0,T],C^2(\R^d)\right)$ 
to the backward Kolmogorov equation~\eqref{eq:PDE}
satisfying the same polynomial growth assumptions at infinity than $f^{\lambda}$ and $g^{\lambda}$
(for instance, see Theorem~5.3 in~\cite{friedman-75}):
\begin{equation}\label{eq:PDE}
\left\{
\begin{array}{l} \displaystyle 
\partial_t u^\lambda + b^\lambda(t,y)\cdot\nabla u^\lambda 
	+ \frac{1}{2} \sigma^\lambda(t,y)\sigma^\lambda(t,y)^T:\nabla^2 u^\lambda = f^{\lambda}(t,y) \ ,
\\
u^\lambda(T,y) = g^{\lambda}(y) \,,
\end{array}
\right.
\end{equation}
where the notation $\nabla u^\lambda$ means $\nabla_y u^\lambda(t,y)$
and $\sigma^\lambda(t,y)\sigma^\lambda(t,y)^T:\nabla^2 u^\lambda$ means 
$\sum_{i,j,k=1}^d\sigma^\lambda_{ik}(t,y)\sigma^\lambda_{jk}(t,y)\partial_{y_i,y_j}^2 u^\lambda(t,y) $.
Using It\^{o} formula for $\left( u^\lambda(t,X_t^{\lambda}) , t \in[0,T] \right)$
with~$u^\lambda$ solution to~\eqref{eq:PDE}, we get the following integral representation of $Z^\lambda$
(see also Appendix~\ref{app:pde_discretization} 
for another link between the SDE~\eqref{eq:pb} and the PDE~\eqref{eq:PDE},
potentially useful to numerics):
\begin{equation}\label{integral_representation}
 g^{\lambda}(X_T^{\lambda})-\int_0^T f^{\lambda}(s,X_s^{\lambda}) ds
= u^\lambda(0,x) + 
\int_0^T \nabla u^\lambda(s,X_s^{\lambda}) \cdot \sigma^{\lambda}(s,X_s^{\lambda}) dB_s \,.
\end{equation}
Note that the left-hand side of~\eqref{integral_representation} is $Z^{\lambda}$, 
and the right-hand side is the sum of a stochatic integral (with zero mean)
plus a scalar $u^\lambda(0,x)$ (thus equal to the expected value $\E{Z^{\lambda}}$ of the left-hand side).
Hence, the optimal control variate also writes:
\begin{equation} \label{eq:ideal_control}
Y^{\lambda}= Z^{\lambda}-\E{Z^{\lambda}} =
\int_0^T \nabla u^\lambda(s,X_s^{\lambda}) \cdot \sigma^{\lambda}(s,X_s^{\lambda}) dB_s \,.
\end{equation}

Of course, the formula~\eqref{eq:ideal_control} is again idealistic because,
most often, numerically solving the PDE~\eqref{eq:PDE} is a very difficult task
(especially in large dimension $d\ge4$).

\subsection{Outline of the algorithms}

Considering either~\eqref{eq:best} or~\eqref{eq:ideal_control},
we propose two algorithms for the efficient online computation of the family 
of parametrized outputs $\{\E{Z^{\lambda}}, \lambda\in\Lambda \}$,
when the parameter $\lambda$ can take any value in a given range $\Lambda$,
using (for each $\lambda\in\Lambda$) 
a control variate built as a linear combination of objects precomputed offline.

More precisely, in \emph{Algorithm~1}, we do the following:
\begin{itemize}
 \item Compute \textit{offline} an accurate approximation $\tilde{Y}^{\lambda}$ of $Y^{\lambda}$
using~\eqref{eq:best},
for a small set of selected parameters $\lambda \in \{\lambda_1, \ldots ,\lambda_I\}\subset\Lambda$ 
(where $I\in\N_{>0}$).
\item For any $\lambda\in\Lambda$, compute \textit{online} 
a control variate for the Monte-Carlo estimation of $\E{Z^{\lambda}}$ 
as a linear combination of $ \{ \tilde Y^{\lambda_i} , i=1,\ldots, I \} \,$:
$$ {\tilde Y}_I^{\lambda} = \sum_{i=1}^I \mu_i^\lambda {\tilde Y}^{\lambda_i} \,. $$
\end{itemize}

And in \emph{Algorithm~2}, we do the following:
\begin{itemize}
 \item Compute \textit{offline} an accurate approximation $\tilde{u}^{\lambda}$ of the solution
${u}^{\lambda}$ to the Kolmogorov backward equation~\eqref{eq:PDE} 
for a small set of selected parameters $\lambda \in \{\lambda_1, \ldots ,\lambda_I\}\subset\Lambda$.
\item For any $\lambda\in\Lambda$, 
compute \textit{online} a control variate for the Monte-Carlo
computation of $\E{Z^{\lambda}}$, 
in view of~\eqref{eq:ideal_control}, 
as a linear combination of 
$\int_0^T \nabla {\tilde u}^{\lambda_i}(s,X_s^{\lambda}) \cdot \sigma^{\lambda}(s,X_s^{\lambda}) dB_s \,$ 
(where $i=1,\ldots, I$):
\begin{equation}
\label{eq:y_approx2}
 {\tilde Y}_I^{\lambda} = 
\sum_{i=1}^I \mu_i^\lambda 
\int_0^T \nabla {\tilde u}^{\lambda_i}(s,X_s^{\lambda}) \cdot \sigma^{\lambda}(s,X_s^{\lambda}) dB_s \,.
\end{equation}
\end{itemize}

For a fixed size $I$ of the reduced-basis, being given a parameter $\lambda$,
both algorithms compute the coefficients $\mu_i^\lambda$, $i=1,\ldots,I$,
with a view to minimizing 
the variance of the random variable $Z^{\lambda}-{\tilde Y}_I^{\lambda}$ 
(in practice, the empirical variance ${\rm Var_M}(Z^{\lambda}-\tilde{Y}_I^{\lambda})$). 

For the moment being, we do not make further precise how we choose the set of parameters 
$\{\lambda_1, \ldots ,\lambda_I\}$ offline. 
This will be done by the same \textit{greedy procedure} for both algorithms, 
and will be the subject of the next section.
Nevertheless,
we would now like to make more precise how we build offline:
\begin{itemize}
 \item[-] in Algorithm~1, approximations $ \{ {\tilde Y}^{\lambda_i} , i=1,\ldots, I
   \} $ for $ \{ Y^{\lambda_i} , i=1,\ldots, I \} $, and
 \item[-] in Algorithm~2, approximations $ \{\nabla{\tilde u}^{\lambda_i} , i=1,\ldots, I
   \} $ for $ \{ \nabla u^{\lambda_i} , i=1,\ldots, I \} $,
\end{itemize}
assuming the parameters $\{ {\lambda_i} , i=1,\ldots, I \}$ have been selected.

For Algorithm~1,
${\tilde Y}^{\lambda_i}$ is built using the fact that it is possible
to compute offline \textit{accurate} Monte-Carlo approximations ${\rm E}_M(Z^{\lambda_i})$ 
of $\E{Z^{\lambda_i}}$ using 
a very large number $M=M_{\rm large}$ of copies of $Z^{\lambda_i}$,
mutually independent and also independant of the copies of $Z^{\lambda}$ used
for the online Monte-Carlo estimation of $\E{Z^{\lambda}}$, $\lambda\neq\lambda_i$
(remember that the amount of offline computations is not meaningful
in the case of a very large number of outputs to be computed online).
The quantities ${\rm E}_{M_{\rm large}}(Z^{\lambda_i})$ 
are just real numbers that can be easily stored in memory at the end of the offline stage 
for re-use online to approximate the  control variate $Y^{\lambda_i}=Z^{\lambda_i}-\E{Z^{\lambda_i}}$ through:
\be \label{eq:y_approx1}
\tilde{Y}^{\lambda_i} = Z^{\lambda_i}-{\rm E}_{M_{\rm large}}(Z^{\lambda_i}) \,.
\ee

For Algorithm~2, we compute approximations $ {\tilde u}^{\lambda_i}$ 
as numerical solutions to the Kolmogorov backward equation~\eqref{eq:PDE}.
For example, in the numerical results of Section~\ref{sec:numerics}, 
the PDE~\eqref{eq:PDE} is solved numerically
with classical deterministic discretization methods
(like finite differences in the calibration problem for instance).

\begin{rem}[Algorithm~2 for stochastic processes with large dimension $d$]
Most deterministic methods to solve a PDE 
(like the finite difference or finite elements methods) remain suitable only for $d\le 3$.
Beyond, one can for example resort to probabilistic discretizations:
namely, a Feynman-Kac representation of the PDE solution,
whose efficiency at effectively reducing the variance has already been shown in~\cite{newton-94}.
We present this alternative probabilistic approximation in Appendix~\ref{app:pde_discretization},
but we will not use it in the present numerical investigation.
\end{rem}

One crucial remark is that for both algorithms,
in the online Monte-Carlo computations,
the Brownian motions which are used to build the control variate (namely $Z^{\lambda_i}$ in~\eqref{eq:y_approx1} for Algorithm~1, and the Brownian motion entering $\tilde Y^\lambda_I$ in~\eqref{eq:y_approx2} for Algorithm~2) 
are the same as those used for~$Z^\lambda$.

Note last that, neglecting the approximation errors
$\tilde{Y}^{\lambda_i}-{Y}^{\lambda_i}$ and $\tilde{u}^{\lambda_i}-{u}^{\lambda_i}$ 
in the reduced-basis elements computed offline, 
a comparison between Algorithms~1 and~2 is possible.
Indeed, remembering the integral representation:
$$ Y^{\lambda_i} 
= \int_0^T \nabla {u}^{\lambda_i}(s,X_s^{\lambda_i})\cdot\sigma^{\lambda_i}(s,X_s^{\lambda_i}) dB_s \,,
$$ 
we see that the reduced-basis approximation of Algorithm~1 has the form:
\begin{equation*}
 {Y}_I^{\lambda} = 
\sum_{i=1}^I \mu_i^\lambda 
\int_0^T \nabla {u}^{\lambda_i}(s,X_s^{\lambda_i}) \cdot \sigma^{\lambda_i}(s,X_s^{\lambda_i}) dB_s \,,
\end{equation*}
while the reduced-basis approximation of Algorithm~2 has the form:
\begin{equation*}
{Y}_I^{\lambda} = 
\sum_{i=1}^I \mu_i^\lambda 
\int_0^T \nabla {u}^{\lambda_i}(s,X_s^{\lambda}) \cdot \sigma^{\lambda}(s,X_s^{\lambda}) dB_s \,.
\end{equation*}
The residual variances $\Var{Y^{\lambda}-{Y}_I^{\lambda}}$ for Algorithms 1 and 2 then respectively read as:
\begin{equation}
\label{eq:error1}
\int_0^T 
 \E{ \left| 
 \nabla u^\lambda \cdot \sigma^{\lambda} (s,X_s^{\lambda})
- \sum_{i=1}^I \mu_i^\lambda 
 \nabla u^{\lambda_i} \cdot \sigma^{\lambda_i}(s,X_s^{\lambda_i}) 
	\right|^2 } ds \,,
\end{equation} 
and:
\begin{equation}
\label{eq:error2}
\int_0^T 
 \E{ \left| 
 \left( \nabla u^\lambda - \sum_{i=1}^I \mu_i^\lambda \nabla u^{\lambda_i} \right)
\cdot \sigma^{\lambda} (s,X_s^{\lambda}) 
 	\right|^2 } ds \,.
\end{equation}
The formulas~\eqref{eq:error1} and~\eqref{eq:error2}
suggest that Algorithm~2 might be more robust than Algorithm~1
with respect to variations of $\lambda$.
This will be illustrated by some numerical results in Section~\ref{sec:numerics}.

\section{Practical variance reduction with approximate control variates}
\label{sec:approximate}

Let us now detail how to select parameters $\{\lambda_i\in\Lambda,i=1,\ldots,I\}$ offline
inside a large \textit{a priori} chosen trial sample $\Lambda_{\rm trial}\subset\Lambda$ of finite size,
and how to effectively compute the coefficients $(\mu_i^\lambda)_{i=1,\ldots,I}$ 
in the linear combinations $\tilde{Y}_I^{\lambda}$
(see Section~\ref{sec:details} for details about practical choices of 
$\Lambda_{\rm trial}\subset\Lambda$). 

\subsection{Algorithm~1}
\label{sec:algo1}

Recall that some control variates $Y^{\lambda}$ are approximated offline with
a computationally expensive Monte-Carlo estimator using ${M_{\rm large}}\gg1$ independent copies of $Z^{\lambda}$:
\begin{equation} \label{eq:H1}
\tilde{Y}^{\lambda}=Z^{\lambda}-E_{M_{\rm large}}(Z^{\lambda})\approx Y^{\lambda}\,,
\end{equation}
for only a few parameters $\{\lambda_i,i=1,\ldots, I\}\subset\Lambda_{\rm trial}$ to be selected.
The approximations $\tilde{Y}^{\lambda_i}$ are then used online
to span a linear approximation space for the set of all control variates $\{ Y^{\lambda} , \lambda\in\Lambda \}$,
thus linear combined as $\tilde{Y}_I^{\lambda}$.
For any $i=1,\ldots,I$, we denote by $\tilde{Y}_i^{\lambda}$ (for any $\lambda\in\Lambda$)
the reduced-basis approximation of ${Y}^{\lambda}$ built as a linear combination of the first $i$ 
selected random variables $\{\tilde{Y}^{\lambda_j},j=1,\ldots, i\}$:
\be \label{eq:Yapproximation1}
\tilde{Y}_i^{\lambda} = \sum_{j=1}^i \mu_j^\lambda \tilde{Y}^{\lambda_j} \approx {Y}^{\lambda}\ ,
\ee
where $(\mu_j^\lambda)_{j=1,\ldots, i}\in\R^i$ is a vector
of coefficients to be computed for each $\lambda$ 
(and each step $i$, but we omit to explicitly denote the dependence of
each entry $\mu_j^\lambda$, $j=1,\ldots,i$, on $i$).
The computation of the coefficients $(\mu_j^\lambda)_{j=1,\ldots, i}$ follows the same procedure offline
(for each step $i=1,\ldots, I-1$) during the reduced-basis construction
as online (when $i=I$): it is based on a variance minimization
principle (see details in Section~\ref{sec:online1}).

With a view to computing $\E{Z^{\lambda}}$ online through computationally \textit{cheap} Monte-Carlo estimations 
${\rm E}_{M_{\rm small}}(Z^{\lambda}-\tilde{Y}_I^{\lambda})$
using only a few ${M_{\rm small}}$ realizations for all $\lambda\in\Lambda$,
we now explain how to select offline a subset $\{\lambda_i,\,i=1,\ldots,I\}\subset\Lambda_{\rm trial}$ 
in order to minimize $\Var{Z^{\lambda}-\tilde{Y}_I^{\lambda}}$ 
(or at least estimators for the corresponding statistical error).

\subsubsection{Offline stage : parameter selection}

\begin{figure} 
\noindent
\fbox{
\begin{minipage}[ht]{\mywidth}
\setlength{\baselineskip}{1.2\baselineskip}
\textbf{Offline}: select parameters $\{\lambda_i\in\Lambda_{\rm trial},i=1,\ldots, I\}$ in
$\Lambda_{\rm trial} \subset \Lambda$ a large finite sample. \\
Selection under stopping criterium: maximal residual variance $\le\varepsilon$.
\begin{tabbing}
\hspace{3mm} \= \hspace{3mm} \= \hspace{3mm} \= \hspace{3mm} \= \kill
\>  \> Let $\lambda_1 \in \Lambda$ be already chosen, \\
\>  \> Compute accurate approximation ${\rm E}_{M_{\rm large}}(Z^{\lambda_1})$ 
of $\E{Z^{\lambda_1}}$. \\
\> \textsf{Greedy procedure:} \\
\>  \> For step $i=1, \ldots, I-1 \ (I>1)$: \\
\>  \>  \> For all $\lambda \in \Lambda_{\rm trial}$, {compute $\tilde{Y}_i^{\lambda}$} as~\eqref{eq:Yapproximation1} 
and (cheap) estimations:\\
\>  \>  \>  \>
$
\epsilon_i(\lambda):={\rm Var}_{M_{\rm small}}\left(Z^{\lambda}-\tilde{Y}_i^{\lambda}\right)\text{ for } 
\Var{Z^{\lambda}-\tilde{Y}_i^{\lambda}} \ .
$
\\
\>  \>  \> Select $\lambda_{i+1} \in 
\underset{\lambda \in\Lambda_{\rm trial}\backslash\{\lambda_j,j=1,\ldots, i\}}{\argmax}
\left\{ \epsilon_i(\lambda) \right\} $.
\\
\>  \>  \> If stopping criterium $\epsilon_i(\lambda_{i+1})\le\varepsilon$, Then Exit \textbf{Offline}.
\\
\>  \>  \> Compute accurate approximation ${\rm E}_{M_{\rm large}}(Z^{\lambda_{i+1}})$ of $\E{Z^{\lambda_{i+1}}}$.
\end{tabbing}
\end{minipage}
}
\caption{Offline stage for Algorithm~1: greedy procedure in metalanguage}
\label{fig:greedy1}
\end{figure}

The parameters $\{\lambda_i,\,i=1,\ldots,I\}$ are selected {\it incrementally}
inside the trial sample $\Lambda_{\rm trial}$
following a \textit{greedy procedure} (see Fig.~\ref{fig:greedy1}).
The incremental search between steps $i$ and $i+1$ reads as follows.
Assume that control variates $\{\tilde{Y}^{\lambda_j},\,j=1,\ldots,i\}$
have already been selected at the step $i$ of the reduced basis construction
(see Remark~\ref{rk:lambda1} for the choice of $\tilde{Y}^{\lambda_1}$).
Then, $\tilde{Y}^{\lambda_{i+1}}$ is chosen following the principle of controlling the 
maximal residual variance inside the trial sample
after the variance reduction using the first $i$ selected random variables:
\begin{equation} \label{eq:criterium1}
\lambda_{i+1}\in\underset{\lambda \in\Lambda_{\rm trial}\backslash\{\lambda_j,j=1,\ldots, i\}}{\mathop{\rm argmax}}
\Var{Z^{\lambda}-\tilde{Y}_i^{\lambda}}\,,
\end{equation}
where the coefficients $(\mu_j^\lambda)_{j=1,\ldots, i}$ 
entering the linear combinations $\tilde{Y}_i^{\lambda}$ in~\eqref{eq:Yapproximation1}
are computed, at each step $i$, 
like for $\tilde{Y}_I^{\lambda}$ in the online stage
(see Section~\ref{sec:online1}).

In practice, the variance in~\eqref{eq:criterium1} is estimated
by an empirical variance~:
$$\Var{Z^{\lambda}-\tilde{Y}_i^{\lambda}}
\simeq{\rm Var}_{\rm M_{small}}(Z^{\lambda}-\tilde{Y}_i^{\lambda})\,.$$
In our numerical experiments, we use the same number ${M_{\rm small}}$ of realizations for the offline computations 
(for all $\lambda \in \Lambda_{\rm trial}$) 
as for the online computations, even though this is not necessary. 
Note that choosing a small number ${M_{\rm small}}$ of realizations 
for the offline computations is advantageous
because the computational cost of the Monte-Carlo estimations in the greedy procedure
is then cheap.
This is useful since $\Lambda_{\rm trial}$ is very large,
and at each step $i$, ${\rm Var}_{\rm M_{small}}(Z^{\lambda}-\tilde{Y}_i^{\lambda})$
has to be computed for all $\lambda\in\Lambda_{\rm trial}$.

Remarkably, after each (offline) step $i$ of the greedy procedure
and for the next online stage when $i=I$,
only a few real numbers should be stored in memory,
namely the collection $\{{\rm E}_{M_{\rm large}}(Z^{\lambda_j}),\,j=1,\ldots,i\}$
along with the corresponding parameters $\{\lambda_j,\,j=1,\ldots,i\}$
for the computation of the approximations~\eqref{eq:H1}.

\begin{rem}\label{rem:criterium}
Another natural criterium for the parameter selection in the greedy procedure 
could be the maximal residual variance \emph{relatively} to the output expectation
\begin{equation} \label{eq:criterium1bis}
\underset{\lambda \in\Lambda_{\rm trial}}{\max}
\frac{\Var{Z^{\lambda}-\tilde{Y}_i^{\lambda}}}{|\E{Z^{\lambda}}|^2}
\simeq 
\underset{\lambda \in\Lambda_{\rm trial}}{\max}
\frac{{\rm Var}_{\rm M_{small}}(Z^{\lambda}-\tilde{Y}_i^{\lambda})}{|{\rm E}_{\rm M_{small}}(Z^{\lambda})|^2}\,.
\end{equation}
This is particularly relevant if the magnitude of the output $\E{Z^{\lambda}}$
is much more sensitive than that of $\Var{Z^{\lambda}}$ to
the variations on $\lambda$.
And it also proved useful for comparison
and discrimination between Algorithms~1 and~2 
in the calibration of
a local parametrized volatility for the Black-Scholes equation 
(see Fig.~\ref{fig:bsdistrib_large}).
\end{rem}

\subsubsection{Online stage : reduced-basis approximation}
\label{sec:online1}

To compute the coefficients $(\mu^\lambda_j)_{j=1,\ldots,i}$ 
in the linear combinations~\eqref{eq:Yapproximation1},
both \textit{online} for any $\lambda\in\Lambda$ when $i=I$
and \textit{offline} for each $\lambda\in\Lambda_{\rm trial}$ and each step $i$
(see greedy procedure above),
we solve a small-dimensional least squares problem 
corresponding to the minimization of (estimators for)
the variance of the random variable $Z^{\lambda}-\tilde{Y}_i^{\lambda}$.

More precisely, in the case $i=I$ (online stage) for instance,
the $I$-dimensional vector $\mu^\lambda = (\mu^\lambda_i)_{1\le i\le I}$ 
is defined, for any $\lambda\in\Lambda$, as the 
unique global minimizer of the following strictly convex
problem of variance minimization:
\begin{equation}\label{eq:minimum1} 
\mu^\lambda = \underset{\mu=(\mu_i)_{1\le i\le I} \in \R^{I}}{\mathop{\rm argmin}} 
\Var{ Z^{\lambda} - \sum_{i=1}^I \mu_i \tilde{Y}^{\lambda_i} 
} \ ,
\end{equation}
or equivalently as the unique solution to the following linear system~:
\begin{equation}\label{eq:linearsystem1}
\sum_{j=1}^I
\Cov{\tilde{Y}^{\lambda_i};\tilde{Y}^{\lambda_j}
} \mu_j^\lambda 
= 
\Cov{\tilde{Y}^{\lambda_i};Z^{\lambda}
}\,,\, \forall i=1, \ldots, I \,.
\end{equation}

Of course, in practice, we use the estimator
(for $X,Y\in L^2_\P(\Omega)$ and $M\in\N_{>0}$)~:
\begin{equation*}
 {\rm Cov_{M}}( X ; Y ) := 
\frac1M \sum_{m=1}^M X_m Y_m
- \left( \frac1M \sum_{m=1}^M X_m \right)
 \left( \frac1M \sum_{m=1}^M Y_m \right)
\end{equation*}
to evaluate the statistical quantities above.
That is, defining a matrix ${\bf C}^{\rm M_{small}}$ with entries 
the following empirical Monte-Carlo estimators ($i,j \in \{1,\ldots, I\}$)~:
$$ \mathbf{C}_{i,j}^{\rm M_{small}}
={\rm Cov_{\rm M_{small}}}\left( \tilde{Y}^{\lambda_i};\tilde{Y}^{\lambda_j}
\right) \,, $$
and a vector ${\bf b}^{\rm M_{small}}$ with entries 
($i\in \{1,\ldots, I\}$)
$ \mathbf{b}_{i}^{\rm M_{small}}
={\rm Cov_{\rm M_{small}}}\left( \tilde{Y}^{\lambda_i};Z^{\lambda}
\right) \,,$
the linear combinations~\eqref{eq:Yapproximation1} are computed 
using as coefficients the Monte-Carlo estimators which are entries of the following vector of $\R^I$:
\begin{equation}\label{eq:mcestimateformu}
\mu^{\rm M_{small}} = \left[{\bf C}^{\rm M_{small}}\right]^{-1} {\bf b}^{\rm M_{small}} \,.
\end{equation}

The cost of one online computation for one parameter $\lambda$
ranges as the computation of $\rm M_{small}$ (independent) realizations
of the random variables $(Z^{\lambda},Y^{\lambda_1},\ldots,Y^{\lambda_I})$,
plus the Monte-Carlo estimators 
${\rm E}_{M_{\rm small}},\,{\rm Cov}_{M_{\rm small}},\,{\rm Var}_{M_{\rm small}}$
and the computation of the solution $\mu^{\rm M_{ small} }$ to the 
(small $I$-dimensional, but full) linear system~\eqref{eq:mcestimateformu}.

In practice, one should be careful when computing~\eqref{eq:mcestimateformu},
because the likely quasi-colinearity of some reduced-basis elements 
often induces ill-conditionning of the matrix ${\bf C}^{\rm M_{small}}$.
Thus the QR or SVD algorithms~\cite{golub-vanloan-96} should be preferred 
to a direct inversion of~\eqref{eq:linearsystem1}
with the Gaussian elimination or the Cholevsky decomposition.
One important remark is that, once the reduced basis is built,
the \textit{same} (small $I$-dimensional) covariance matrix
${\bf C}^{\rm M_{small}}$
has to be inverted for \textit{all} $\lambda\in\Lambda$,
as soon as the same Brownian paths are used for each online evaluation. 
And the latter condition is easily satisfied in practice,
simply by resetting the seed of the random number generator to the same value 
for each new online evaluation (that is for each new $\lambda\in\Lambda$).

\begin{rem}[Final output approximations and bounds]
It is a classical result that,
taking first the limit ${\rm M_{large}}\to\infty$ then ${\rm M_{small}}\to\infty$,
$ \mu^{\rm M_{small}} \overset{\P-a.s.}{\underset{{\rm M_{small},M_{large}}\to\infty}{\longrightarrow}} \mu^\lambda $.
So, the variance 
is indeed (asymptotically) reduced
to the minimum $\Var{Z^{\lambda}-{Y}^{\lambda}_I}$ in~\eqref{eq:minimum1},
obtained with the optimal linear combination ${Y}^{\lambda}_I$ of
selected control variates ${Y}^{\lambda_i}$ (without approximation).
In addition, using Slutsky theorem twice successively for Monte-Carlo
estimators of the coefficient vector $\mu^\lambda$ and of the variance $\Var{Z^{\lambda}-{Y}^{\lambda}_I}$,
it also holds a computable version of the Central Limit Theorem,
which is similar to~\eqref{eq:clt}
except that it uses Monte-Carlo estimations of 
$Z^{\lambda}-{\tilde Y}^{\lambda}_I$
instead of $Z^{\lambda}$ to compute the confidence intervals
(and with successive limits ${\rm M_{large}}\to\infty$, 
${\rm M_{small}}\to\infty$).
So our output approximations now read for all $\lambda\in\Lambda$:
$$ \E{Z^\lambda} \simeq {\rm E}_{M_{\rm small}}\left(
Z^{\lambda} - \sum_{i=1}^I \mu_i^{\rm M_{small}} \tilde{Y}^{\lambda_i} 
\right)\,, $$
and asymptotic probabilistic error bounds are given by the confidence intervals~\eqref{eq:clt}.
\end{rem}

\subsection{Algorithm~2}
\label{sec:algo2}

In Algorithm~2,
approximations $\nabla\tilde{u}^{\lambda_i}$ 
of the gradients $\nabla u^{\lambda_i}$ of the solutions $u^{\lambda_i}$ 
to the backward Kolmogorov equation~\eqref{eq:PDE} are computed offline
for only a few parameters $\{\lambda_i,i=1,\ldots, I\}\subset\Lambda_{\rm trial}$ to be selected.
In comparison with Algorithm~1, 
approximations $(\nabla\tilde{u}^{\lambda_i})_{i=1,\ldots,I}$ are now used online
to span a linear approximation space for $\{\nabla{u}^{\lambda}\,,\ \lambda\in\Lambda\}$.
At step $i$ of the greedy procedure ($i=1,\ldots,I$),
the reduced-basis approximations $\tilde{Y}_i^{\lambda}$ 
for the control variates $Y^{\lambda}$ 
read (for all $\lambda\in\Lambda$):
\begin{align} 
\label{eq:Yapproximation2}
\tilde{Y}_i^{\lambda}
&=\sum_{j=1}^i\mu_j^\lambda \tilde{Y}^{\lambda_j}_\lambda  \approx {Y}^{\lambda}\,,
\\
\label{eq:Yapproximation2rb}
\tilde{Y}^{\lambda_j}_\lambda 
&=  \int_0^T\nabla\tilde{u}^{\lambda_j}(s,X_s^{\lambda})\cdot\sigma^{\lambda}(s,X_s^{\lambda}) dB_s \,.
\end{align}
where $(\mu_j^\lambda)_{j=1,\ldots,i}$ are coefficients to be computed 
for each $\lambda$ (again, the dependence of $\mu_j^\lambda$ on the step $i$ is implicit).
Again, the point is to explain, first, 
how to select parameters $\{\lambda_i,\,i=1,\ldots,I\}\subset\Lambda_{\rm trial}$ in the offline stage,
and second, how to compute the coefficients $(\mu_j^\lambda)_{j=1,\ldots,i}$
in each of the $i$-dimensional linear combinations $\tilde{Y}_i^{\lambda}$.
Similarly to Algorithm~1, 
the parameters $\{\lambda_i,\,i=1,\ldots,I\}\subset\Lambda_{\rm trial}$ are selected offline 
following the greedy procedure,
and, for any $i=1,\ldots,I$, the coefficients $(\mu_j^\lambda)_{j=1,\ldots,i}$ 
in the linear combinations offline and online
are computed, following the same principle of minimizing the variance,
by solving a least squares problem.

\subsubsection{Offline stage : parameter selection}

\begin{figure} 
\noindent
\fbox{
\begin{minipage}[ht]{\mywidth}
\setlength{\baselineskip}{1.2\baselineskip}
\textbf{Offline}: select parameters $\{\lambda_i\in\Lambda_{\rm trial},i=1,\ldots, I\}$ in
$\Lambda_{\rm trial} \subset \Lambda$ a large finite sample. \\
Selection under stopping criterium: maximal residual variance $\le\varepsilon$.
\begin{tabbing}
\hspace{3mm} \= \hspace{3mm} \= \hspace{3mm} \= \hspace{3mm} \= \kill
\>  \> Let $\lambda_1 \in \Lambda$ be already chosen, \\
\>  \> Compute approximation $\nabla \tilde{u}^{\lambda_1}$ of $\nabla u^{\lambda_1}$. \\
\> \textsf{Greedy procedure:} \\
\>  \> For step $i=1, \ldots, I-1 \ (I>1)$: \\
\>  \>  \> For all $\lambda \in \Lambda_{\rm trial}$, 
{compute $\tilde{Y}_i^{\lambda}$} 
as~\eqref{eq:Yapproximation2} and estimations:\\
\>  \>  \>  \>
$
\epsilon_i(\lambda):={\rm Var}_{M_{\rm small}}\left(Z^{\lambda}-\tilde{Y}_i^{\lambda}\right)
\text{ for } \Var{Z^{\lambda}-\tilde{Y}_i^{\lambda}} \ .
$
\\
\>  \>  \> Select $\lambda_{i+1} \in 
\underset{\lambda \in\Lambda_{\rm trial}\backslash\{\lambda_j,j=1,\ldots, i\}}{\argmax}
\left\{ \epsilon_i(\lambda) \right\} $.
\\
\>  \>  \> If stopping criterium $\epsilon_i(\lambda_{i+1})\le\varepsilon$, Then Exit \textbf{Offline}.
\\
\>  \>  \> Compute approximation
$\nabla \tilde{u}^{\lambda_{i+1}}$ of $\nabla u^{\lambda_{i+1}}$.
\end{tabbing}
\end{minipage}
}
\caption{Offline stage for Algorithm~2: greedy procedure in metalanguage}
\label{fig:greedy2}
\end{figure}

The selection of parameters $\{\lambda_j,\,j=1,\ldots,i\}$
from a trial sample $\Lambda_{\rm trial}$
follows a greedy procedure like in Algorithm~1 (see Fig.~\ref{fig:greedy2}).
In comparison with Algorithm~1, after $i$ (offline) steps of the greedy procedure 
($1 \le i \le I-1$) and online ($i=I$),
note that discretizations of functions $(t,y)\to\nabla\tilde{u}^{\lambda_j}(t,y)$,
$j=1,\ldots,i+1$, 
are stored in memory to compute
the stochastic integrals~\eqref{eq:Yapproximation2}, 
which is possibly a huge amount of data.

\subsubsection{Online stage : reduced-basis approximation}

Like in Algorithm~1, the coefficients $(\mu^\lambda_j)_{j=1,\ldots,i}$ 
in the linear combination~\eqref{eq:Yapproximation2}  are computed similarly 
\textit{online} (and then $i=I$) for any $\lambda\in\Lambda$ 
and \textit{offline} (when $1\le i \le I-1$) for each $\lambda\in\Lambda_{\rm trial}$ 
as minimizers of -- a Monte Carlo discretization of -- the 
least squares problem:
\begin{equation}\label{eq:minimum2}
\min_{ \mu \in \R^{I} } 
\Var{ Z^{\lambda} - \sum_{i=1}^I \mu_i \tilde{Y}^{\lambda_i}_\lambda } \ ,
\end{equation}
where we recall that $\tilde{Y}^{\lambda_i}_\lambda$
are defined by~\eqref{eq:Yapproximation2rb}.
Note that contrary to the reduced-basis elements $\tilde{Y}^{\lambda_i}$ in Algorithm~1,
the elements $\tilde{Y}^{\lambda_i}_\lambda$ in Algorithm~2 have to be 
recomputed for \textit{each} queried parameter value $\lambda\in\Lambda$.

Again, in practice,
the unique solution  $(\mu^\lambda_j)_{j=1,\ldots,i}$ to the variational problem~\eqref{eq:minimum2} 
is equivalently
the unique solution to the following linear system:
\begin{equation}\label{eq:linearsystem2}
\sum_{j=1}^I
\Cov{\tilde{Y}^{\lambda_i}_\lambda;\tilde{Y}^{\lambda_j}_\lambda} \mu_j^\lambda = 
\Cov{\tilde{Y}^{\lambda_i}_\lambda;Z^{\lambda}}\,,\, \forall i=1,\ldots,I \,,
\end{equation}
and is in fact computed as the unique solution to the discrete minimization problem:
\begin{equation}
\label{eq:minimum2discrete}
\mu^{\rm M_{small}} = \left[{\bf C}^{\rm M_{small}}\right]^{-1} {\bf b}^{\rm M_{small}}
\,,
\end{equation}
with 
$ \mathbf{C}_{i,j}^{\rm M_{small}}
={\rm Cov_{\rm M_{small}}}\left( \tilde{Y}^{\lambda_i}_\lambda;
\tilde{Y}^{\lambda_j}_\lambda\right) $
and 
$ \mathbf{b}_{i}^{\rm M_{small}}
={\rm Cov_{\rm M_{small}}}\left( \tilde{Y}^{\lambda_i}_\lambda;Z^{\lambda}
\right) \,.$

The cost of one computation online for one parameter $\lambda$ 
is more expensive than that in Algorithm~1,
and ranges as the computation of $M_{\rm small}$ independent realizations of $Z^{\lambda}$,
\emph{plus} the computation of $I$ (discrete approximations of) 
the stochastic integrals~\eqref{eq:Yapproximation2rb},
plus the Monte-Carlo estimators and the solution $\mu^{\rm M_{small}}$ to the 
(small $I$-dimensional, but full) linear system~\eqref{eq:minimum2discrete}.
In comparison to Algorithm~1,
notice that the (discrete) covariance matrix ${\bf C}^{\rm M_{small}}$ to be inverted 
depends on $\lambda$, and thus cannot be treated offline once for all:
it has to be recomputed for each $\lambda\in\Lambda$.

\subsection{General remarks about reduced-basis approaches}

The success of our two reduced-basis approaches clearly depends on 
the variations of $Z^{\lambda}$ with $\lambda\in\Lambda$.
Unfortunately, we do not have yet 
a precise 
understanding of this,
similarly to the PDE case~\cite{patera-rozza-07}.
Our reduced-basis approaches
have only been investigated numerically in relevant cases for application 
(see Section~\ref{sec:numerics}).
So we now provide some theoretical ground 
only for the {\it a priori} existence of a reduced basis,
like in the PDE case~\cite{maday-patera-turinici-02},
with tips for a practical use of the greedy selection procedure
based on our numerical experience.
Of course, it remains to show that
the greedy procedure actually selects a good reduced basis.

\subsubsection{{\it A priori} existence of a reduced basis}
\label{sec:apriori}

Following the analyses~\cite{maday-patera-turinici-02,patera-rozza-07} 
for parametrized PDEs,
we can prove the {\it a priori} existence of a reduced basis
for some particular collections of parametrized control variates,
under very restrictive assumptions on the structure of the parametrization.
\begin{prop}\label{prop:apriori_rb}
Assume there exist collections of uncorrelated
(parameter-independent) random variables with zero mean 
$Y_j \in L^2_\P(\Omega)$, $1\le j\le J$,
and of positive $C^\infty(\R)$ 
functions $g_j$, $1\le j\le J$, such that 
\begin{equation} \label{eq:Yassumption}
Y^\lambda = \sum_{j=1}^J g_j(\lambda) \, Y_j \,,\ \forall \lambda \in \Lambda \,,
\end{equation}
and there exists a constant $C>0$ such that, 
for all parameter ranges $\Lambda = [\lambda_{\min},\lambda_{\max}] \subset \R$,
there exists a $C^\infty$ diffeomorphism $\tau_\Lambda$ defined on $\Lambda$ satisfying:
\begin{equation}\label{eq:gassumption}
\sup_{1\le j\le J} \sup_{\tilde{\lambda}\in\tau_{\Lambda}(\Lambda)}
 (g_j\circ\tau_\Lambda^{-1})^{(M)}(\tilde{\lambda}) \le M! C^M \,, 
\text{ for all $M$-derivatives of $g_j\circ\tau_\Lambda^{-1}$.} 
\end{equation}
Then, for all parameter ranges $\Lambda = [\lambda_{\min},\lambda_{\max}] \subset \R$, 
there exist constants $c_1,c_2>0$ independent of $\Lambda$ and $J$ such that,
for all $N\in\mathbb{N}_{>0}$, 
$ N\ge N_0:= 1+c_1\left(\tau_{\Lambda}(\lambda_{\max})-\tau_\Lambda(\lambda_{\min})\right) $,
there exist 
$N$ distinct parameter values 
$ \lambda_n^N\in\Lambda\,,\ n=1,\ldots,N\,,\ (\text{with }\lambda_n^N\le\lambda_{n+1}^N) \,,$ 
sastisfying, with $\mathcal{Y}_N=\Span{Y^{\lambda_n^N},n=1,\ldots,N}$: 
\begin{equation} \label{eq:rbcv} 
\inf_{Y_N \in \mathcal{Y}_N} \Var{Z^\lambda-Y_N} \le 
e^{-\frac{c_2}{N_0-1}(N-1)} \, 
\Var{Z^\lambda} 
\,,\ {\forall\lambda\in\Lambda} \,.
\end{equation}
\end{prop}
One can always write $Y^\lambda$ like~\eqref{eq:Yassumption}
with uncorrelated random variables (using a Gram-Schmidt procedure)
and with positive coefficients 
(at least on a range $\Lambda$ where they do not vanish).
But the assumption~\eqref{eq:gassumption} is much more restrictive.
The mapping $\tau_\Lambda$ for the parameter,
which depends on the functions $g_j$, $j=1,\ldots,J$,
indeed tells us 
how the convergence depends on variations in the size of the parameter range~$\Lambda$.
See~\cite{maday-patera-turinici-02,patera-rozza-07}
for an example of such functions $g_j$ and $\tau_\Lambda$,
and Appendix~\ref{app:apriori_rb} for a short proof inspired 
from~\cite{maday-patera-turinici-02,patera-rozza-07}.

The Proposition~\ref{prop:apriori_rb}
may cover a few interesting cases of application
for the {\it a priori} existence theory.
One example where the assumption~\eqref{eq:Yassumption} 
hold is the following.
Consider an output $Z^\lambda=g(X_T^\lambda)$ with $g$ a polynomial function,
and~:
\begin{equation}\label{eq:pbou}
X_t^{\lambda} = x + \int_0^t b^\lambda(s) X_s^{\lambda} \, ds 
	+ \int_0^t \sigma^{\lambda}(s) dB_s \,.
\end{equation}
The optimal control variate $Y^\lambda$
in such a case writes in the form~\eqref{eq:Yassumption}
(to see this, one can first explicitly compute
the reiterated (or multiple) It\^o integrals 
in the polynomial expression of $g(X_T^\lambda)$
with Hermite polynomials~\cite{karatzas-shreve-91}).
Then,~\eqref{eq:gassumption} may hold 
provided $b^\lambda$ and $\sigma^\lambda$ 
are smooth functions of $\lambda\in\Lambda$
(again, see~\cite{maday-patera-turinici-02,patera-rozza-07}
for functions $g_j$ satisfying~\eqref{eq:gassumption}).
But quite often, 
the reduced bases selected in practice by the {\it greedy} procedure
are much better than $\mathcal{Y}_N$
(see~\cite{patera-rozza-07} for comparisons
when $\lambda$ is scalar).

\subsubsection{Requirements for efficient practical greedy selections}
\label{sec:details}

A comprehensive study would clearly need hypotheses 
about the regularity of $Y^\lambda$ as a function of $\lambda$ 
and about the discretization $\Lambda_{\rm trial}$ of $\Lambda$
to show that the greedy procedure actually selects good reduced bases.
We do not have precise results yet,
but we would nevertheless like to provide the reader
with conjectured requirements for the greedy procedure to work
and help him as a potential user of our method.

Ideally, one would use the greedy selection procedure 
directly on $\{Y^\lambda,\lambda\in\Lambda\}$ for Algorithm~1 
and on $\{\nabla u^\lambda,\lambda\in\Lambda\}$ for Algorithm~2.
But in pratice, one has to resort to approximations only,
$\{{\tilde Y}^\lambda,\lambda\in\Lambda\}$ for Algorithm~1
and $\{{\nabla \tilde u}^\lambda,\lambda\in\Lambda\}$ for Algorithm~2.
So,
following requirements on discretizations of parametrized PDEs in the classical reduced-basis method~\cite{patera-rozza-07},
the stability of the reduced basis selected by the greedy procedure for parametrized control variates intuitively requires:
\begin{itemize}
 \item[(H1)] 
\textit{
For any required accuracy $\varepsilon>0$, 
we assume the existence of approximations,
$\tilde{Y}^{\lambda}$ for $Y^{\lambda}$ in Algorithm~1
(resp. $\tilde{u}^{\lambda}$ for $u^{\lambda}$ in Algorithm~2),
such that the $L^2$-approximation error is uniformly bounded on $\Lambda$:
}
\end{itemize}
\begin{gather*}
\forall\lambda\in\Lambda\,, \E{|\tilde{Y}^{\lambda}-Y^{\lambda}|^2}\le\varepsilon \,,
\\
\left(\textit{resp. }
\int_0^T\E{|\nabla\tilde{u}^{\lambda}-\nabla u^{\lambda}|^2(X^\lambda_t)}dt
\le \varepsilon
\text{ or }
\|\nabla\tilde{u}^{\lambda}-\nabla{u}^{\lambda}\|^2_{L^2}\le\varepsilon\right)\,. 
\end{gather*}
Moreover, in practice, one can only manipulate finite nested samples of parameter $\Lambda_{\rm trial}$ instead of the full range $\Lambda$.
So some representativity assumption about $\Lambda_{\rm trial}$ 
is also intuitively required
for the greedy selection procedure to work on $\Lambda$:
\begin{itemize}
 \item[(H2)] 
\textit{
For any required accuracy $\varepsilon>0$,
we assume the existence of a sufficiently \textit{representative} finite discrete subset 
$\Lambda_{\rm trial}\subset\Lambda$ of parameters such that
reduced bases built from $\Lambda_{\rm trial}$ are still good enough for $\Lambda$.
}
\end{itemize}
Refering to Section~\ref{sec:apriori}, {\it good enough} reduced bases
should satisfy exponential convergence like~\eqref{eq:rbcv},
with slowly deteriorating capabilities in terms of approximation
when the size of the parameter range grows.
Now, in absence of more precise result, 
intuitition has been necessary so far to choose good discretizations.
The numerical results of Section~\ref{sec:numerics}
have been obtained with $\rm M_{large} = 100\:M_{small}$ in Algorithm~1,
and with a trial sample $\Lambda_{\rm trial}$ of $100$ parameter values
randomly chosen (with uniform distribution) in $\Lambda$.

In absence of theory for the greedy procedure,
one could also think of using another parameter selection procedure 
in the offline stage.
The interest of the greedy procedure is that it is cheap 
while effective in practice.
In comparison, 
another natural reduced basis would be defined by the first $I$ leading eigenvectors
from the Principal Components Analysis (PCA) 
of the very large covariance matrix with entries 
$\Cov{Y^{\lambda_i};Y^{\lambda_j}}_{(\lambda_i,\lambda_j)
\in\Lambda_{\rm trial}\times\Lambda_{\rm trial}}$.
The latter (known as the Proper Orthogonal Decomposition method)
may yield similar variance reduction 
for most parameter values $\lambda\in\Lambda$~\cite{patera-rozza-07},
but would certainly require more computations during the offline stage.

\begin{rem}\label{rk:lambda1}
The choice of the first selected parameter $\lambda_1$ 
has not been precised yet.
It is observed that most often, this choice does not impact the quality of the variance reduction.
But to be more precise, we choose $\lambda_1\in\Lambda_{\rm small\ trial}$
such that $Z^{\lambda_1}$ has maximal variance 
in a small initial sample $\Lambda_{\rm small\ trial}\subset\Lambda$, 
for instance.
\end{rem}

\section{Worked examples and numerical tests}
\label{sec:numerics}

The efficiency of our reduced-basis strategies for parametrized problems
is now investigated numerically for two problems relevant to some applications.

\begin{rem}[High-dimensional parameter]\label{rem:highdim}
Although the maximal dimension in the parameter treated here is two,
one can reasonably hope for our reduced-basis approach to 
remain feasible with moderately high-dimensions in the parameter range $\Lambda$,
say twenty.
Indeed, a careful mapping of a multi-dimensional parameter range
may allow for an efficient sampling $\Lambda_{\rm trial}$
that makes a greedy procedure tractable and
next yields a good reduced basis for $\Lambda$,
as it was shown for the classical reduced-basis method with
parametrized PDEs~\cite{sen-08,boyaval-lebris-maday-nguyen-patera-08}.
\end{rem}

\subsection{Scalar process with constant drift and parametrized diffusion}
\label{sec:bs}

\subsubsection{Calibration of the Black--Scholes model with local volatility}
\label{sec:bs1}

One typical computational problem in finance is the valuation of
an option depending on a risky asset with value $S_t$ at time $t\in[0,T]$.
In the following we consider Vanilla European Call options 
with payoff $\phi(S_T;K)=\max(S_T-K,0)$,
$K$ being the exercise price (or strike) of the option at time $t=T$.
By the no arbitrage principle for a portfolio mixing 
the risky asset of value $S_t$ with a riskless asset of interest rate $r(t)$,
the price (as a function of time) is a martingale given by a conditional expectation:
\begin{equation}\label{eq:bsform}
e^{-\int_t^T r(s) ds} \E{\phi(S_T)|\mathcal{F}_t}
\end{equation}
where, in the Black-Scholes model with local volatility,
$S_t = S_t^\lambda$ is a stochastic process solving the Black-Scholes equation:
\begin{equation}\label{eq:bseq}
dS_t^\lambda = S_t^\lambda\left( r(t) \, dt + \sigma^\lambda(t,S_t^\lambda)\:d B_t\right)
\qquad S_{t=0}^\lambda = S_0 \,,
\end{equation}
and $(\mathcal{F}_t)$ is the natural filtration for the standard Brownian motion $(B_t)$.
For this model to be predictive
the parameter $\lambda$ in the (local) volatility $\sigma^\lambda$
needs to be calibrated against observed data.

Calibration, like many numerical optimization procedures, defines a typical many-query context,
where one has to compute many times the price~\eqref{eq:bsform} of the option
for a large number of parameter values
until, for some optimal parameter value $\lambda$,
a test of adequation with statistical data $P_{({K},{\bar t_l})}$ 
observed on the market at times $\bar t_l\in[0,T]\,,\ l=0,\ldots,\bar L$ is satisfied.
For instance, a common way to proceed is to minimize in $\lambda$ the quadratic quantity:
$$\mathcal{J}(\lambda)= \sum_{l=0}^{\bar L} 
\left| e^{-\int_{\bar t_l}^T r(s) \, ds} 
\E{ \phi(S_{T}^\lambda;K) | \mathcal{F}_{\bar t_l} } -  P_{({K},{\bar t_l})}\right|^2\,,$$
most often regularized with some Tychonoff functional,
using optimization algorithms like descent methods which indeed require 
many evaluations of the functional $\mathcal{J}(\lambda)$ for various~$\lambda$.
One could even consider the couple $(K,T)$ as additional parameters 
to optimize the contract, but we do not consider such an extension here.

Note that the reduced-basis method for parameterized 
PDEs~\cite{machiels-maday-patera-01,maday-patera-turinici-02,patera-rozza-07}
has recently proved very efficient at treating a 
similar calibration problem~\cite{pironneau-09}.
Our approach is different since we consider a probabilistic pricing numerical method.

In the following numerical results, we solve~\eqref{eq:bsform} for many parameter values
assuming that the interest rate $r$ is a fixed given constant
and the local volatility $\sigma^\lambda$ has ``hyperbolic'' parametrization~\eqref{eq:param}
(used by practitionners in finance):
\begin{equation} \label{eq:param}
\sigma^\lambda(t,S)= (\Gamma+1) \left( \frac{1}{C(0,S_0)}+\frac{\Gamma}{C(t,S)} \right)^{-1}
\end{equation}
where $C(t,S)=\frac12\left(\sqrt{C_A(t,S)^2+C_{\rm min}^2}+C_A(t,S)\right)$ with:
$$
C_A(t,S)=a+\frac12\sqrt{(b-c)^2\log^2\left(\frac{S}{\alpha S_0 e^{r t}}\right) + 4 a^2 d^2}
+\frac12 (b+c)\log\left(\frac{S}{\alpha S_0 e^{r t}}\right)\,.
$$
The local volatility $\sigma^\lambda$ is thus parametrized with 
a $7$-dimensional parameter $\lambda=(a,b,c,d,\alpha,\Gamma,C_{\rm min})$.

Our reduced-basis approach aims at building a vector space 
in order to approximate the family of random variables:
$$\left\{ Y^\lambda:=e^{-rT}\max(S_T^\lambda-K,0)- e^{-rT}\E{\max(S_T^\lambda-K,0)},\lambda\in\Lambda\right\}\,,$$
which are optimal control variates for the computation of the expectation of $e^{-rT}\max(S_T^\lambda-K,0)$.
In Algorithm~2, 
we also use the fact that
\begin{equation} \label{eq:bsY}
Y^\lambda 
= \int_0^T \partial_S u^\lambda(t,S_t^\lambda) \sigma^\lambda(t,S_t^\lambda) S_t^\lambda d B_t \,, 
\end{equation}
where the function $u^\lambda(t,S)$ solves 
for $(t,S)\in[0,T)\times(0,\infty)$:
\begin{equation} \label{eq:kolmobs}
\partial_t u^\lambda(t,S) + r S \partial_S u^\lambda(t,S) 
+ \frac{\sigma^\lambda(t,S)^2 S^2}{2} \partial_{SS} u^\lambda(t,S) = 0 \,, 
\end{equation}
with final condition $u^\lambda(T,S)=e^{-rT}\max(S-K,0)$.
Note the absence of boundary condition at $S=0$ because 
the advection and diffusion terms are zero at $S=0$.
The backward Kolmogorov equation~\eqref{eq:kolmobs}
is \textit{numerically} solve using 
finite differences~\cite{achdou-pironneau-05}.
More precisely, after a change of variable $u^\lambda(t,S)=e^{-rt}C^\lambda(t,S)$,
equation~\eqref{eq:bsY} rewrites:
\begin{equation}\label{eq:vcbs}
Y^\lambda = \int_0^T 
e^{-rt} \partial_S C^\lambda(t,S_t^\lambda) \sigma^\lambda(t,S_t^\lambda) S_t^\lambda d B_t \,,
\end{equation}
where  $C^\lambda(t,S)$ solves the classical Black-Scholes PDE:
\begin{equation}\label{eq:calleq}
\partial_t C^\lambda(t,S) - r C^\lambda(t,S) + r S \partial_S C^\lambda(t,S) + \frac{\sigma^\lambda(t,S)^2 S^2}{2} \partial_{SS} C^\lambda(t,S) = 0 \,,
\end{equation}
with the final condition $C^\lambda(T,S)=\max(S-K,0)$.
In the case of a low-dimensional variable~$S_t$ (like one-dimensional here),
one can use 
a finite differences method of order 2 (with Crank-Nicholson discretization in time) 
to compute approximations 
$\tilde{C}^\lambda_{l,j} \simeq C^\lambda(t_l,x_j)$, $l=0,\ldots,L$, $j=0,\ldots,J$
on a grid for the truncated domain $[0,T]\times[0,3K]\subset[0,T]\times[0,\infty)$,
with $L=100$ steps in time and $J=300$ steps in space of constant sizes
(and with Dirichlet boundary condition 
$\tilde{C}^\lambda_{l,J+1}=(3-e^{-r(T-t_l)})K\,,\forall l=0,\ldots,N$
at the truncated boundary).
An approximation $\tilde{C}^\lambda(t,S)$ of $C^\lambda(t,S)$ 
at any $(t,S)\in[0,T]\times[0,3K]$ 
is readily reconstructed as a linear interpolation
on tiles $(t,S)\in[t_l,t_{l+1}]\times[S_j,S_{j+1}]$.

\subsubsection{Numerical results}

\begin{figure}[htbp]
\centering
\includegraphics[trim = 0mm 0mm 0mm 0mm, clip, scale=0.8]{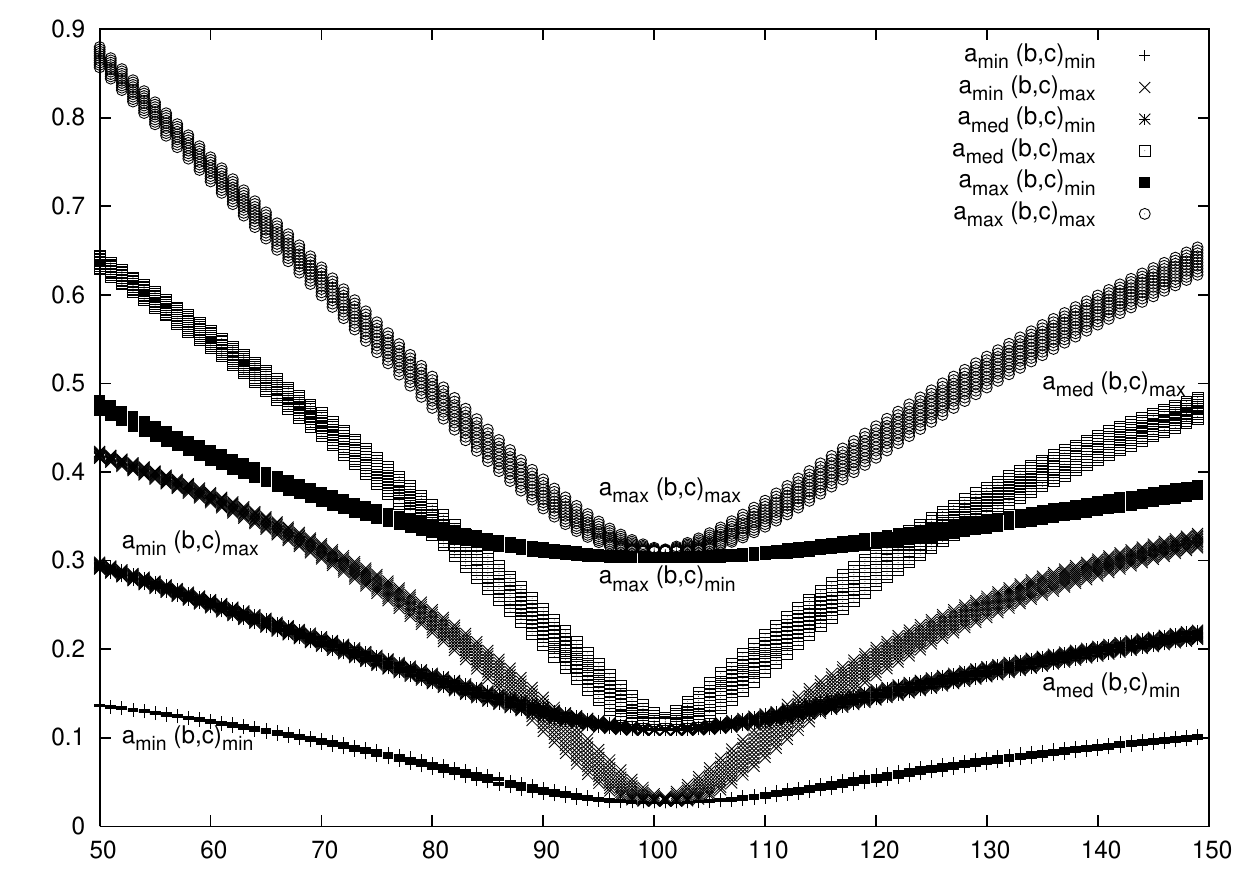}
\caption{Variations of the ``hyperbolic'' local volatility function $\sigma^\lambda(t,S)$
with respect to $S\in[50,150]$.
Six families of curves are shown (as time $t$ evolves in $[0,1]$)
for extremal and mid- values of the parameter $(a,b=c)$ in $[-.05,.15]\times\{b=c\in[.5,1.5]\}$:
$(\min(a),\min(b=c))$, $(\min(a),\max(b=c))$,
$({\rm med}(a):=.5\min(a)+.5\max(a),\min(b=c))$, $({\rm med}(a):=.5\min(a)+.5\max(a),\max(b=c))$,
$(\max(a),\min(b=c))$, $(\max(a),\max(b=c))$.
Each family of curves shows the time variations of $S\to\sigma^\lambda(t,S)$
for $t\in\{.1\times k | k=0,\ldots, 10\}$).
\label{fig:vol}
}
\end{figure}

\begin{figure}[htbp]
\includegraphics[trim = 8mm 0mm 20mm 10mm, clip, scale=.37]{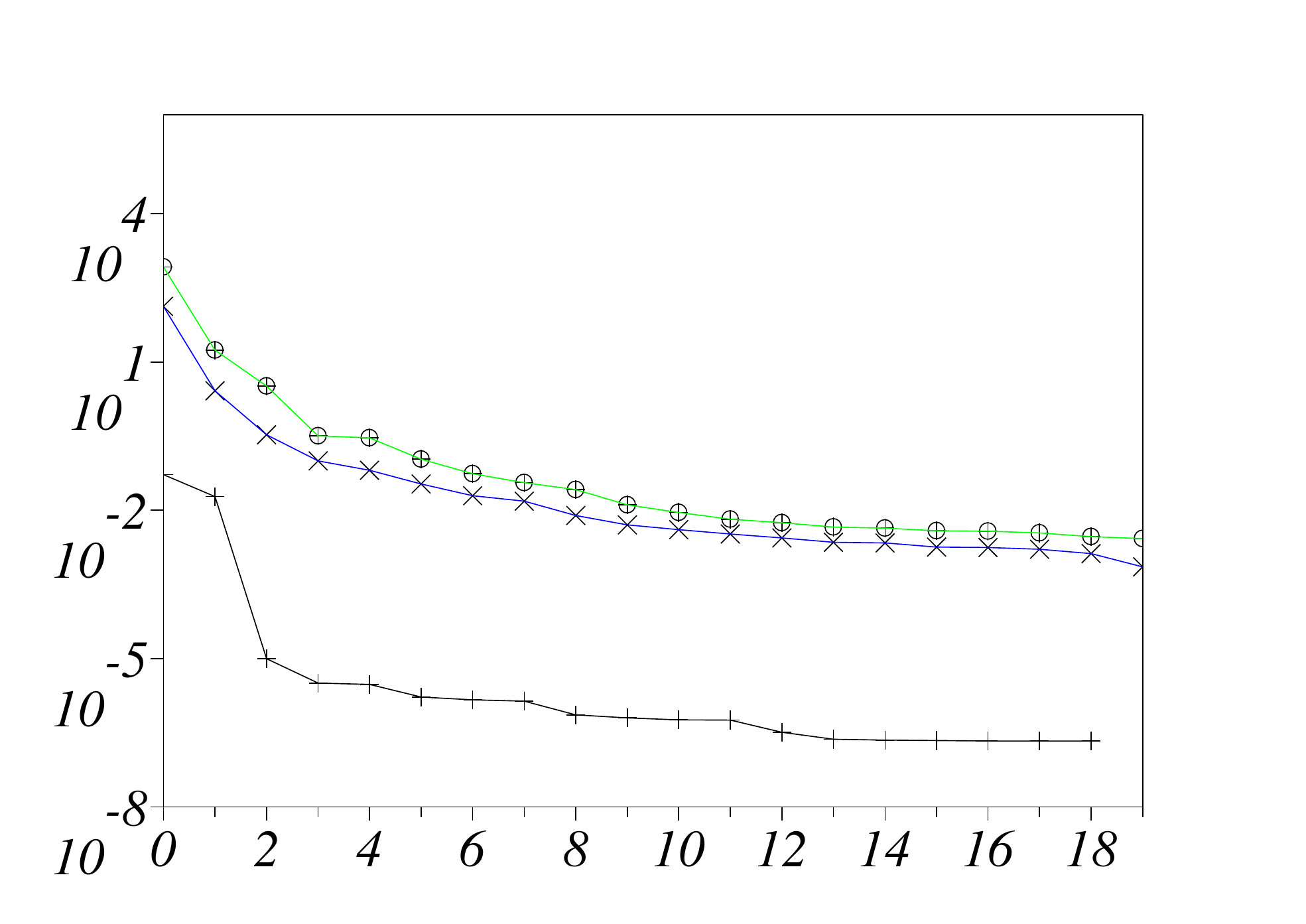}
\includegraphics[trim = 8mm 0mm 20mm 10mm, clip, scale=.37]{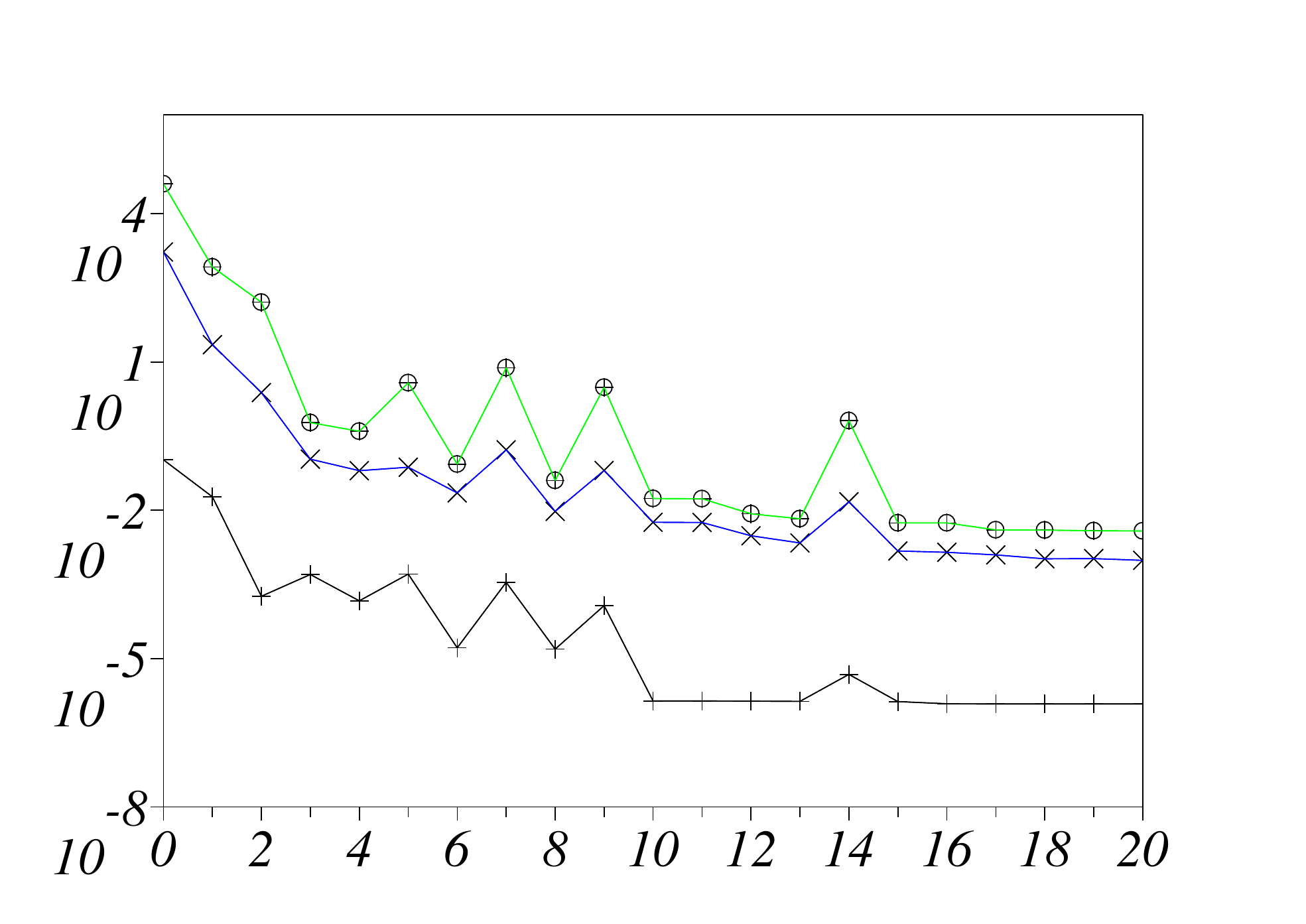}
\caption{Algorithm~1 for Black--Scholes model with local ``hyperbolic'' volatility: 
Minimum $+$, mean $\times$ and maximum $\circ$ of the absolute variance~\eqref{eq:absoluterv} 
in samples of parameters
(left: offline sample $\Lambda_{\rm trial}\setminus\{\lambda_i,i=1,\ldots,I\}$; right: online sample $\Lambda_{\rm test}$) 
with respect to the size $I$ of the reduced basis.
\label{fig:bsdistrib1}
}
\includegraphics[trim = 8mm 0mm 20mm 10mm, clip, scale=.37]{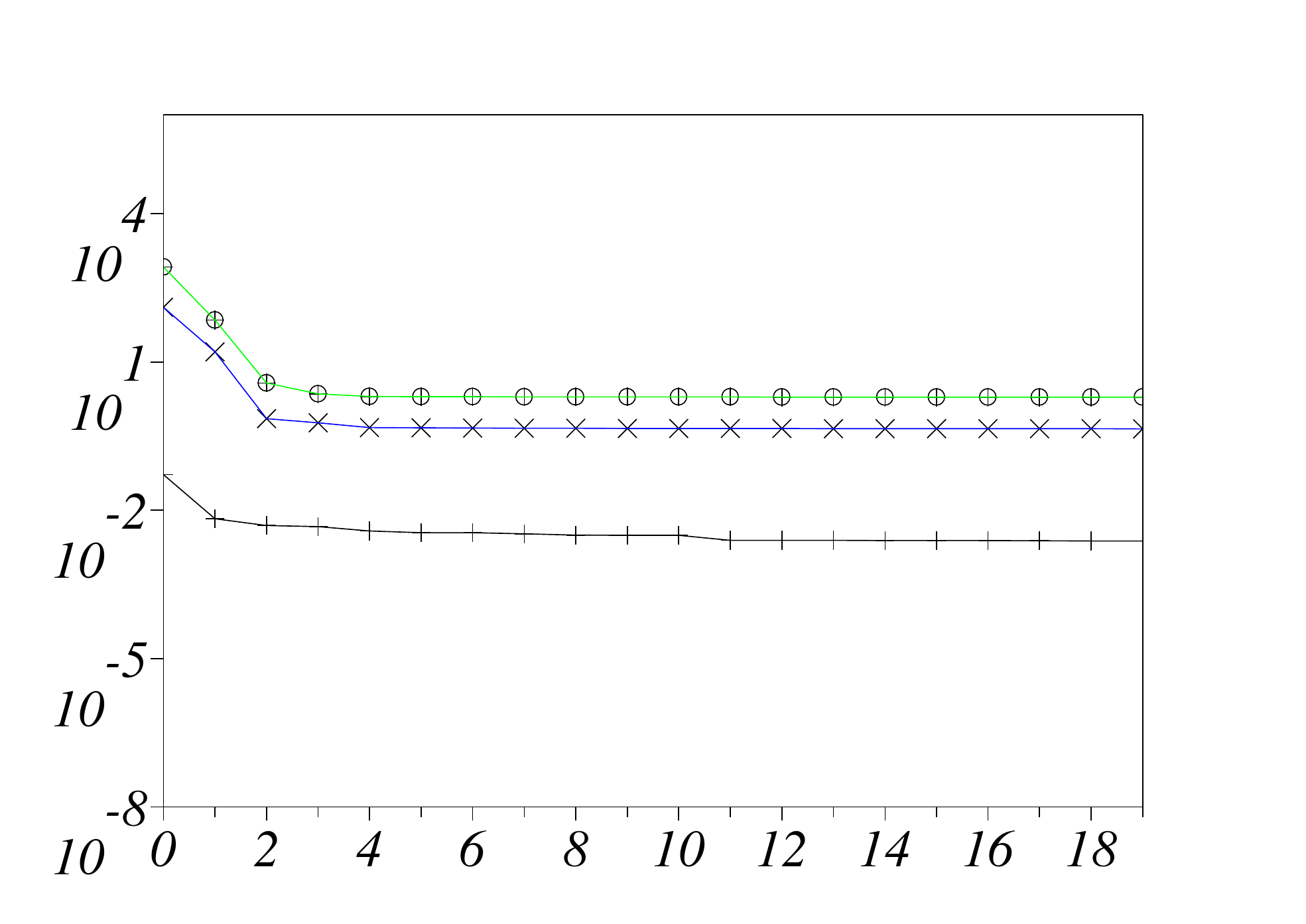}
\includegraphics[trim = 8mm 0mm 20mm 10mm, clip, scale=.37]{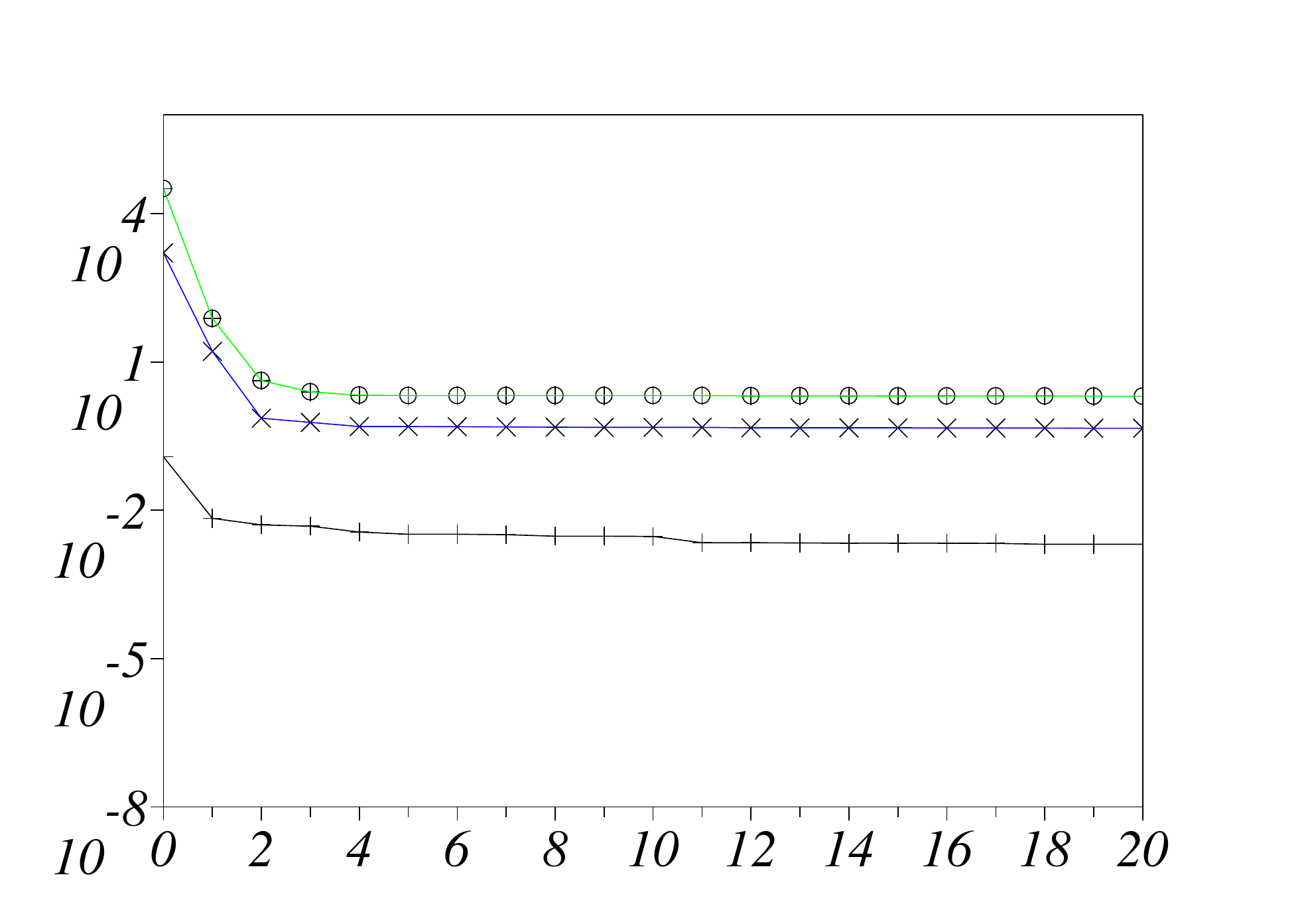}
\caption{Algorithm~2 for Black--Scholes model with local ``hyperbolic'' volatility: 
Minimum $+$, mean $\times$ and maximum $\circ$ of the absolute variance~\eqref{eq:absoluterv} 
in samples of parameters
(left: offline sample $\Lambda_{\rm trial}\setminus\{\lambda_i,i=1,\ldots,I\}$; right: online sample $\Lambda_{\rm test}$) 
with respect to the size $I$ of the reduced basis.
\label{fig:bsdistrib1_2}
}
\includegraphics[trim = 8mm 0mm 20mm 10mm, clip, scale=.37]{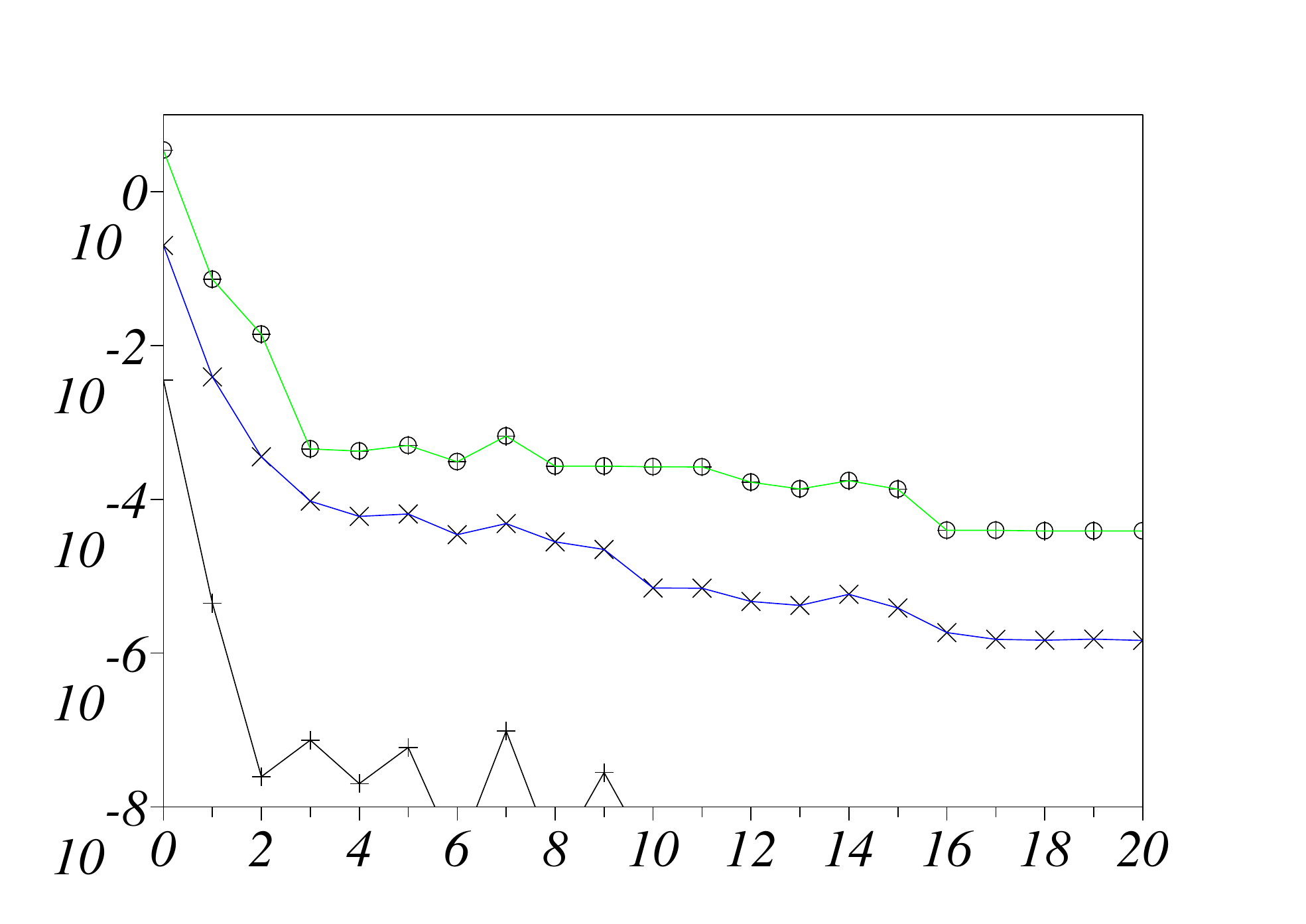}
\includegraphics[trim = 8mm 0mm 20mm 10mm, clip, scale=.37]{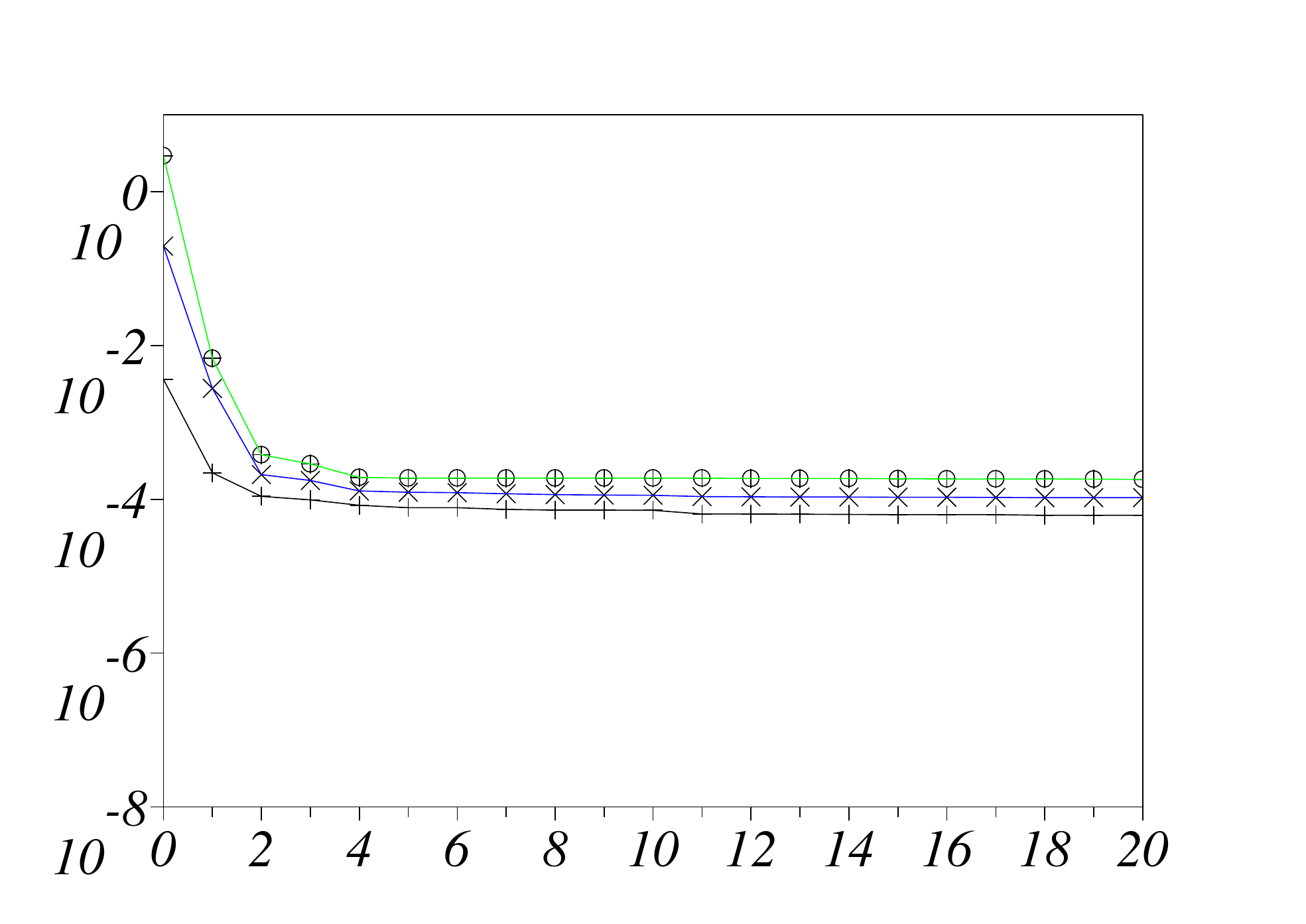}
\caption{Algorithm~1 (left) and 2 (right) for Black--Scholes model with local ``hyperbolic'' volatility: 
Minimum $+$, mean $\times$ and maximum $\circ$ of the relative variance~\eqref{eq:relativerv} 
in a sample test (online) $\Lambda_{\rm test}$ of parameters
with respect to the size $I$ of the reduced basis.
\label{fig:bsdistrib2}
}
\end{figure}

The Euler-Maruyama scheme with $N=10^2$ time steps of constant size $\dt=\frac{T}{N}=10^{-2}$
is used to compute one realization of a pay-off $\max(\tilde{S}_{N}^{\lambda}-K,0)$,
for a strike $K=100$ at final time $t_N=T=1$ 
when the initial price is $\tilde{S}_{0}^{\lambda}=90$
and the interest rate $r=0.04$.
Then, (a large number of) expectations $\E{\max(\tilde{S}_{N}^{\lambda}-K,0)}$ are approximated
through Monte-Carlo evaluations 
$E_{\rm M_{\rm small}}\left(\max(\tilde{S}_{N}^{\lambda}-K,0)\right)$
with ${\rm M_{\rm small}}=10^3$ realizations,
when the local volatility parameter $\lambda=(a,b,c,d,\alpha,\Gamma,C_{\rm min})$ 
assumes many values in the two-dimensional range
$\Lambda=[-.05,.15]\times\{b=c\in[.5,1.5]\}\times\{1.\}\times\{1.1\}\times\{5\}\times\{.05\}$
(variations of the function $\sigma^\lambda(t,S)$ with $\lambda$ are shown in Fig.~\ref{fig:vol}).
We build reduced bases of different sizes $I=1,\ldots,20$ 
from the same sample $\Lambda_{\rm trial}$
of size $|\Lambda_{\rm trial}|=100$,
either with Algorithm~1 (Fig.~\ref{fig:bsdistrib1} and~\ref{fig:bsdistrib2})
using approximate control variates
computed with 
${\rm M_{large}}=100\,{\rm M_{small}}$ evaluations~:
$$ \tilde{Y}^{\lambda}_I = \sum_{i=1}^I \mu_i^{\rm M_{small}} 
\tilde{Y}^{\lambda_i}
= \sum_{i=1}^I\mu_i^{\rm M_{small}} 
\left( \max(\tilde{S}_{N}^{\lambda_i}-K,0)-
E_{\rm M_{\rm large}}\left(\max(\tilde{S}_{N}^{\lambda_i}-K,0)\right) \right)\,,
$$
or with Algorithm~2 (Fig.~\ref{fig:bsdistrib1_2} and~\ref{fig:bsdistrib2})
using approximate control variates: 
$$ \tilde{Y}^{\lambda}_I 
=  \sum_{i=1}^I \mu_i^{\rm M_{small}} 
\left( \sum_{n=0}^{N-1}
e^{-r t_n} \partial_S \tilde{C}^{\lambda_i} (t_n,\tilde{S}_{n}^\lambda)  \sigma^\lambda(t_n,\tilde{S}_{n}^\lambda) \sqrt{|t_{n+1}-t_{n}|}\, G_n \right) $$
computed as first-order discretizations of the It\^o stochastic integral~\eqref{eq:vcbs}
using the finite-difference approximation of 
the solution to the backward Kolmogorov equation.
We always start the greedy selection procedure
by choosing $\lambda_1$ such that $ \tilde{Y}^{\lambda_1}$ 
has the maximal correlation with other members in $\Lambda_{\rm small\ trial}$,
a small prior sample of $10$ parameter values chosen randomly 
with uniform law in $\Lambda$, see Remark~\ref{rk:lambda1}.

We show in Fig.~\ref{fig:bsdistrib1} and~\ref{fig:bsdistrib1_2}
the absolute variance after variance reduction~:
\begin{equation} \label{eq:absoluterv}
{\rm Var_{M_{small}}}\left(
\max(\tilde{S}_{N}^{\lambda}-K,0)-\tilde{Y}^{\lambda}_{I}
\right)\,,
\end{equation}
and in Fig.~\ref{fig:bsdistrib2} 
the relative variance after variance reduction~:
\begin{equation} \label{eq:relativerv}
\frac{ {\rm Var_{\rm M_{small}}}\left(
\max(\tilde{S}_{N}^{\lambda}-K,0)-\tilde{Y}^{\lambda}_{I}
\right)}{
{\rm E_{\rm M_{small}}}\left(
\max(\tilde{S}_{N}^{\lambda}-K,0)-\tilde{Y}^{\lambda}_{I}
\right)^2
}\,.
\end{equation}
In each figure,
the maximum, the minimum and the mean of one of the two residual variance above is shown,
either within the offline sample deprived of the selected parameter values
$\Lambda_{\rm trial}\setminus\{\lambda_i,\,i=1,\ldots,I\}$,
or within an online uniformly distributed sample test $\Lambda_{\rm test}\subset\Lambda$
of size $|\Lambda_{\rm test}|=10\,|\Lambda_{\rm trial}|$.

It seems that Algorithm~1 slighlty outperfoms Algorithm~2 
with a sufficiently large reduced basis,
comparing the (online) decrease rates 
for
either the relative variance or the absolute variance.
Yet, one should also notice that, with very small-dimensional reduced basis,
the Algorithm~2 yields very rapidly good variance reduction.
Comparing the decrease rates of the variance in offline and online samples
tells us how good was the (randomly uniformly distributed here) choice of 
$\Lambda_{\rm trial}$.
The Algorithm~2 seems more robust than the Algorithm~1
for reproducing (``extrapolating'')
offline results from a sample $\Lambda_{\rm trial}$ in the whole range $\Lambda$.
So, comparing the first results for Algorithms~1 and~2,
it is not clear which algorithm performs the best variance reduction
for a given size of the reduced basis.

\begin{figure}[htbp]
\includegraphics[trim = 8mm 0mm 20mm 10mm, clip, scale=.37]{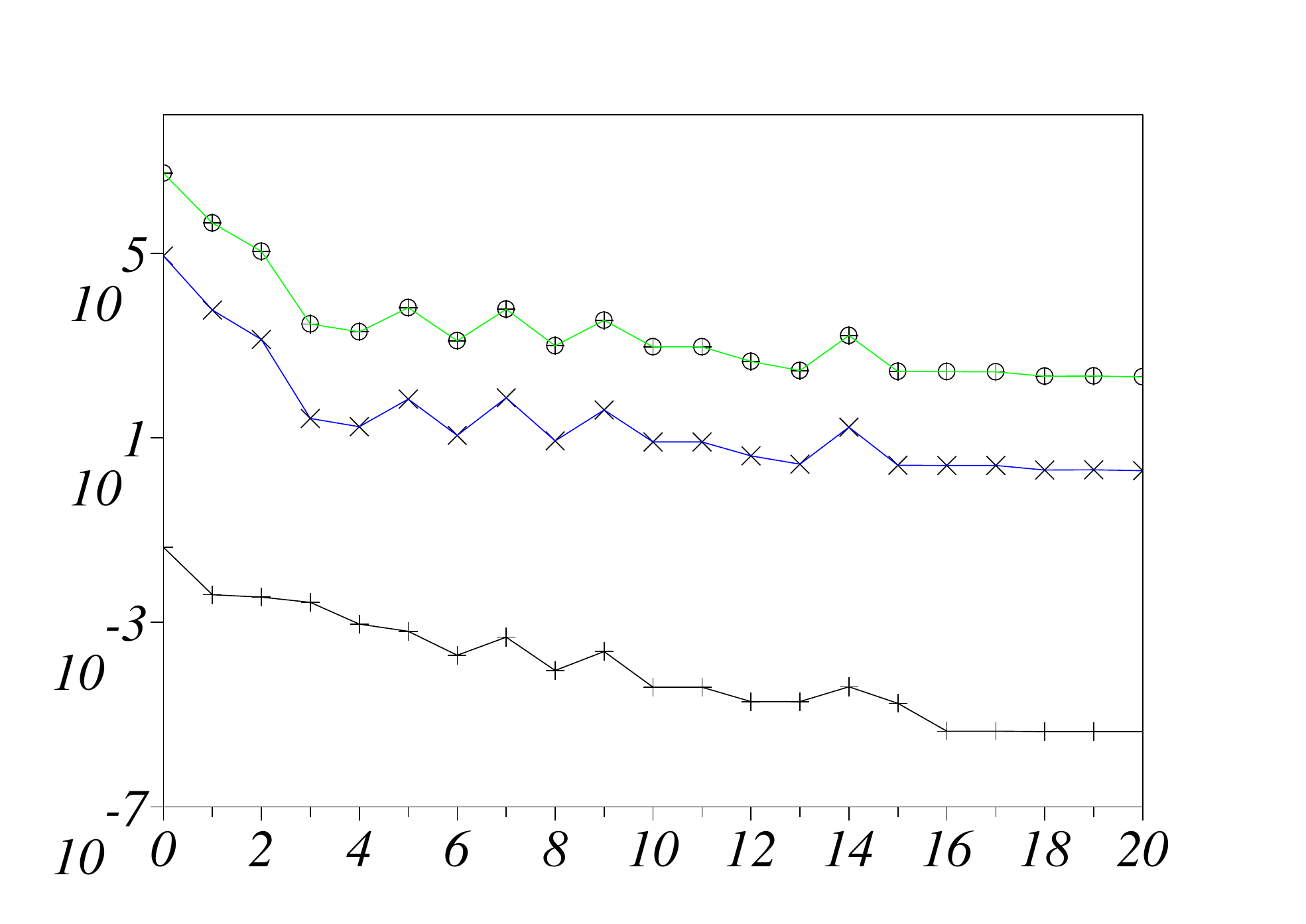}
\includegraphics[trim = 8mm 0mm 20mm 10mm, clip, scale=.37]{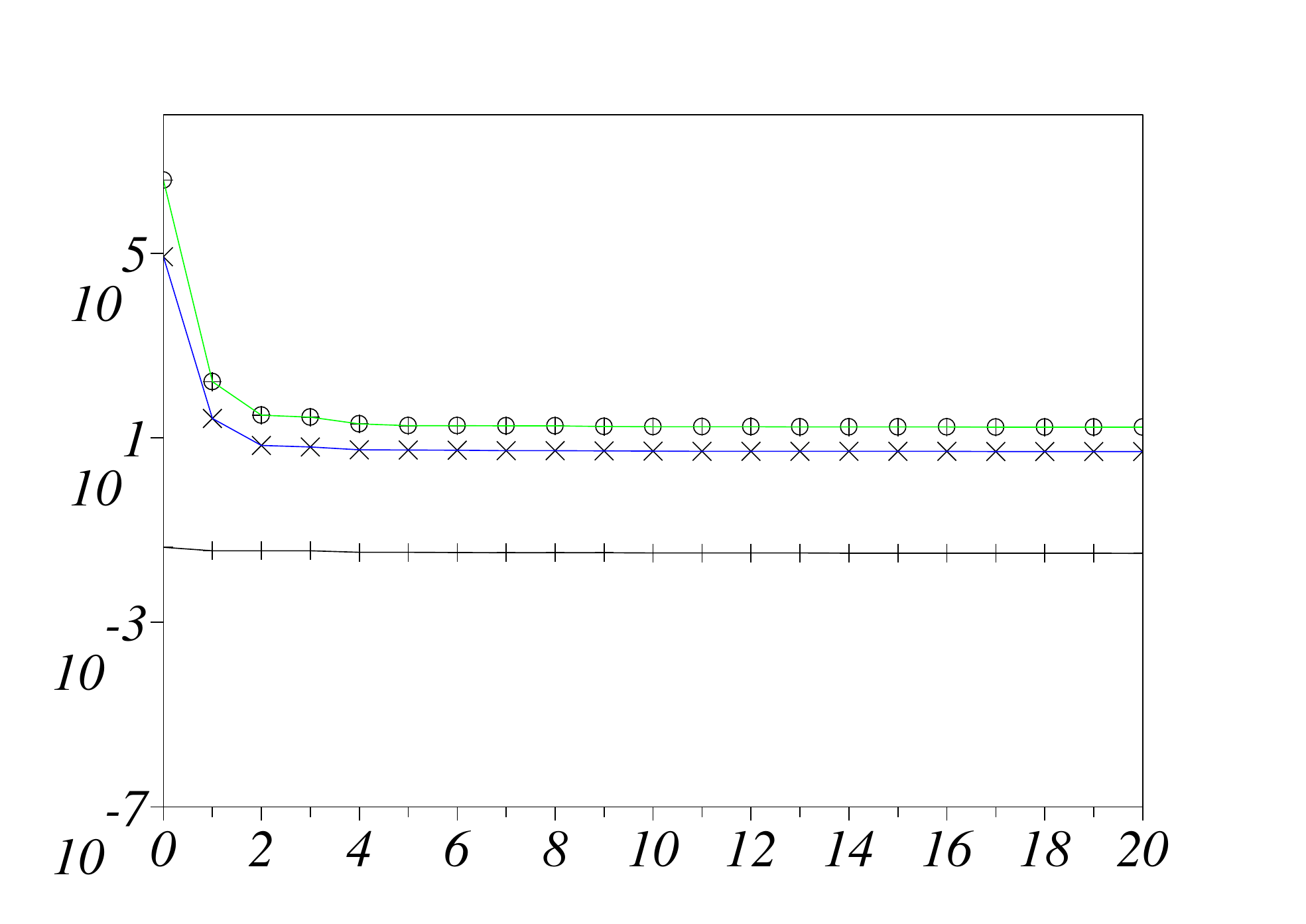}
\\
\includegraphics[trim = 8mm 0mm 20mm 10mm, clip, scale=.37]{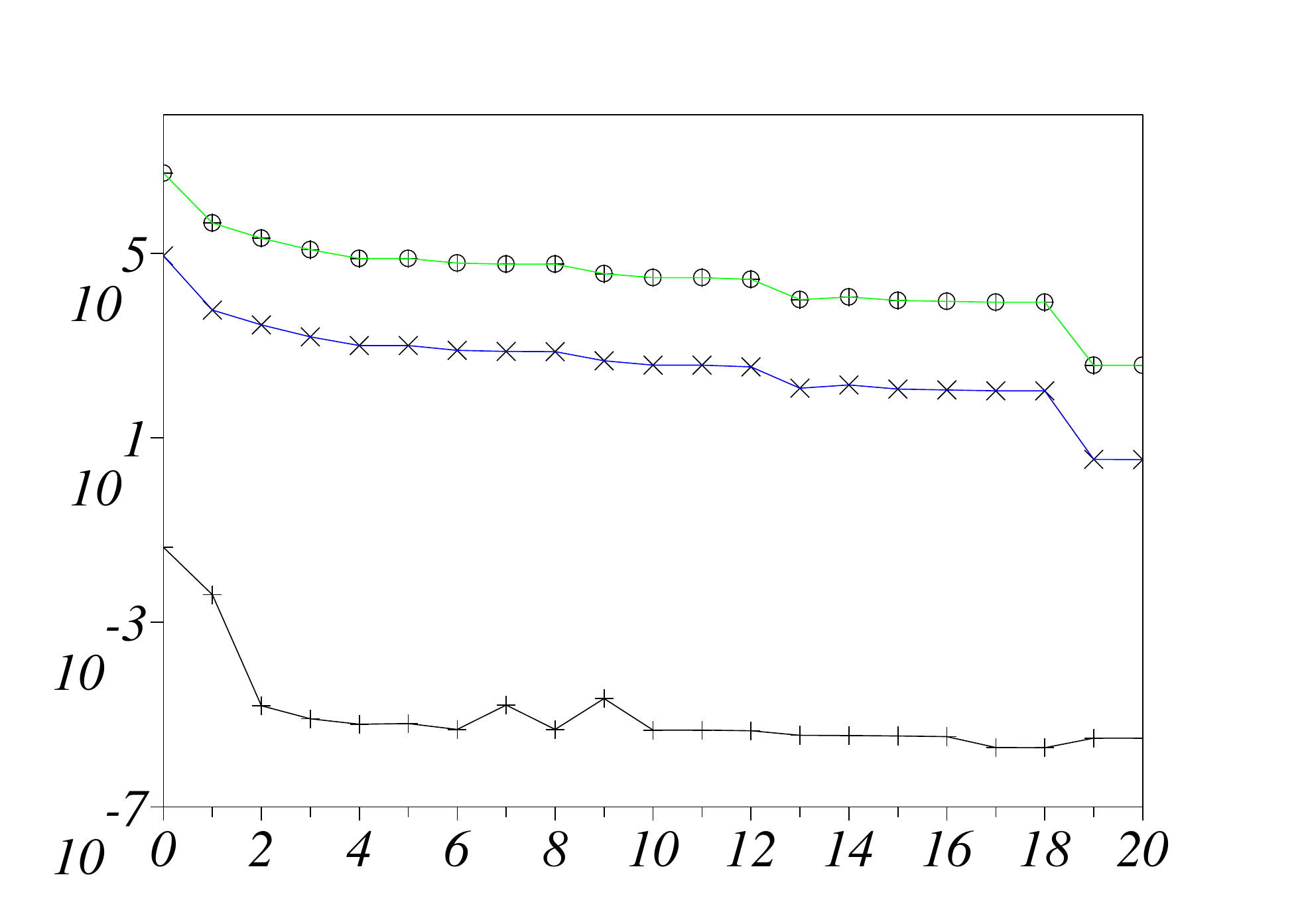}
\includegraphics[trim = 8mm 0mm 20mm 10mm, clip, scale=.37]{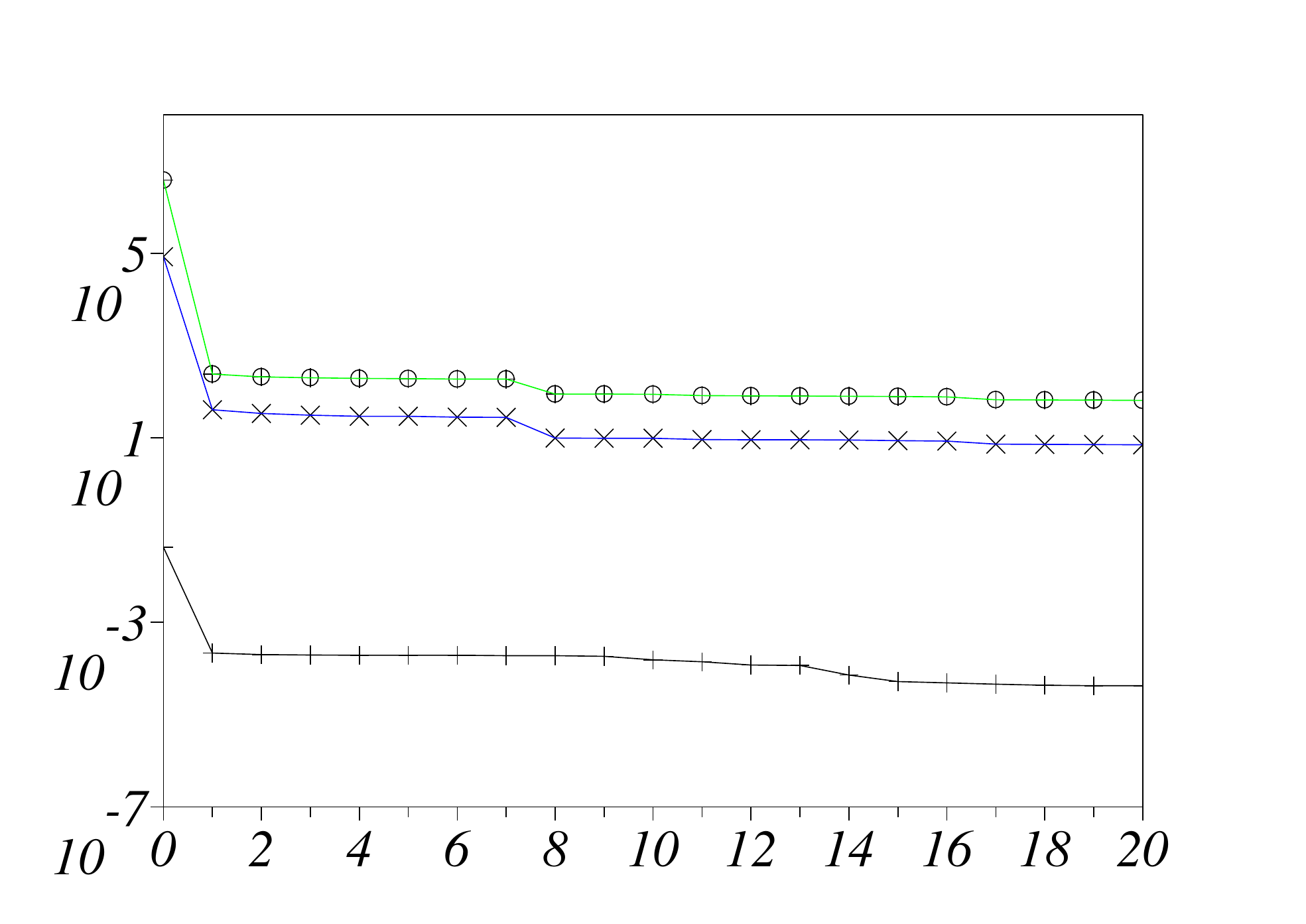}
\caption{Algorithm~1 (left) and 2 (right) for Black--Scholes model with local ``hyperbolic'' volatility: 
Minimum $+$, mean $\times$ and maximum $\circ$ of the 
(online) absolute variance~\eqref{eq:absoluterv} 
in a sample test $\Lambda_{\rm test wide}$ of parameters
with respect to the size $I$ of the reduced basis.
Greedy selection with absolute variance~\eqref{eq:absoluterv} (top) and 
relative variance~\eqref{eq:relativerv} (bottom).
}
\label{fig:bsdistrib_large}
\end{figure}

\begin{figure}[htbp]
\includegraphics[trim = 8mm 0mm 20mm 10mm, clip, scale=.37]{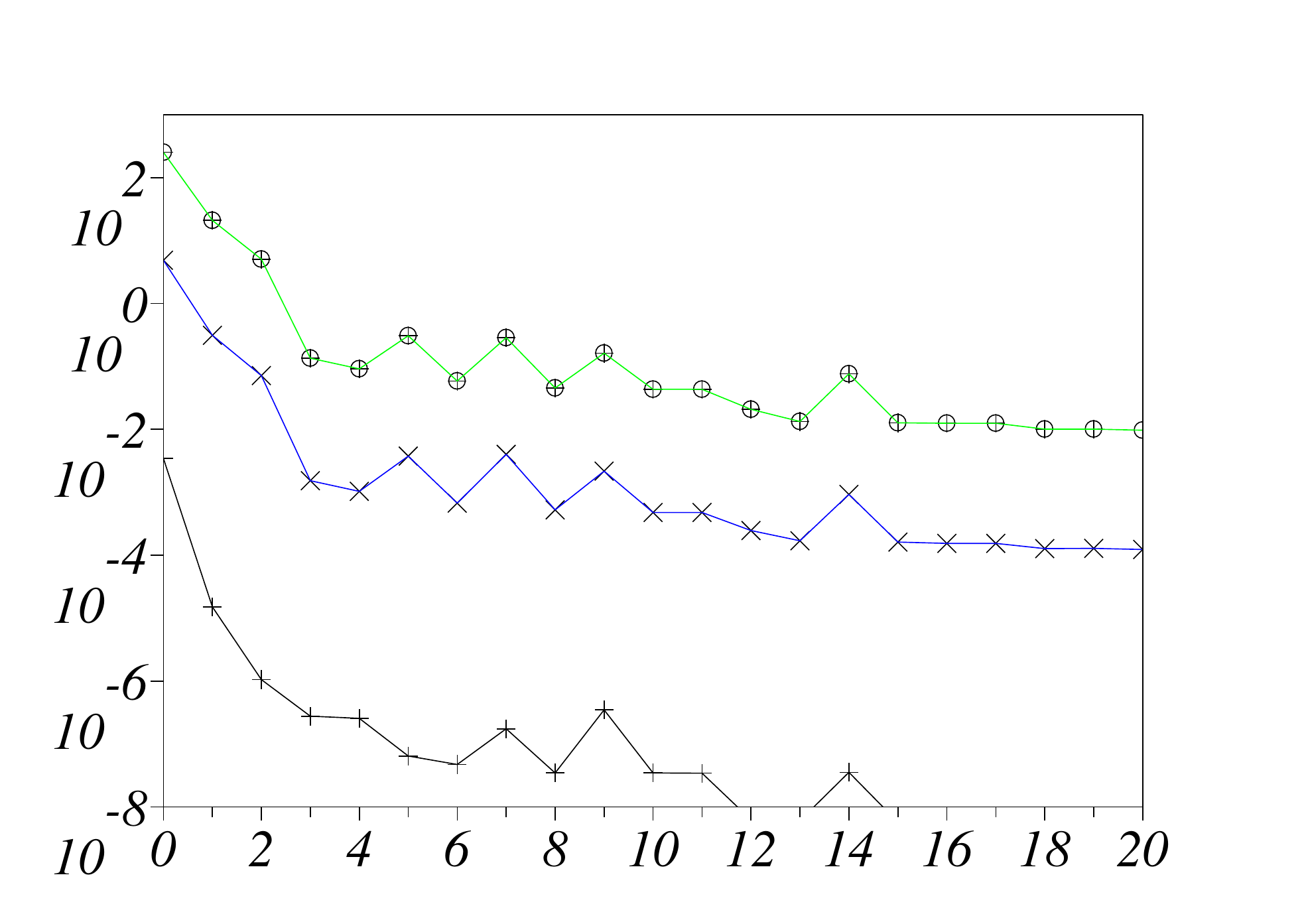}
\includegraphics[trim = 8mm 0mm 20mm 10mm, clip, scale=.37]{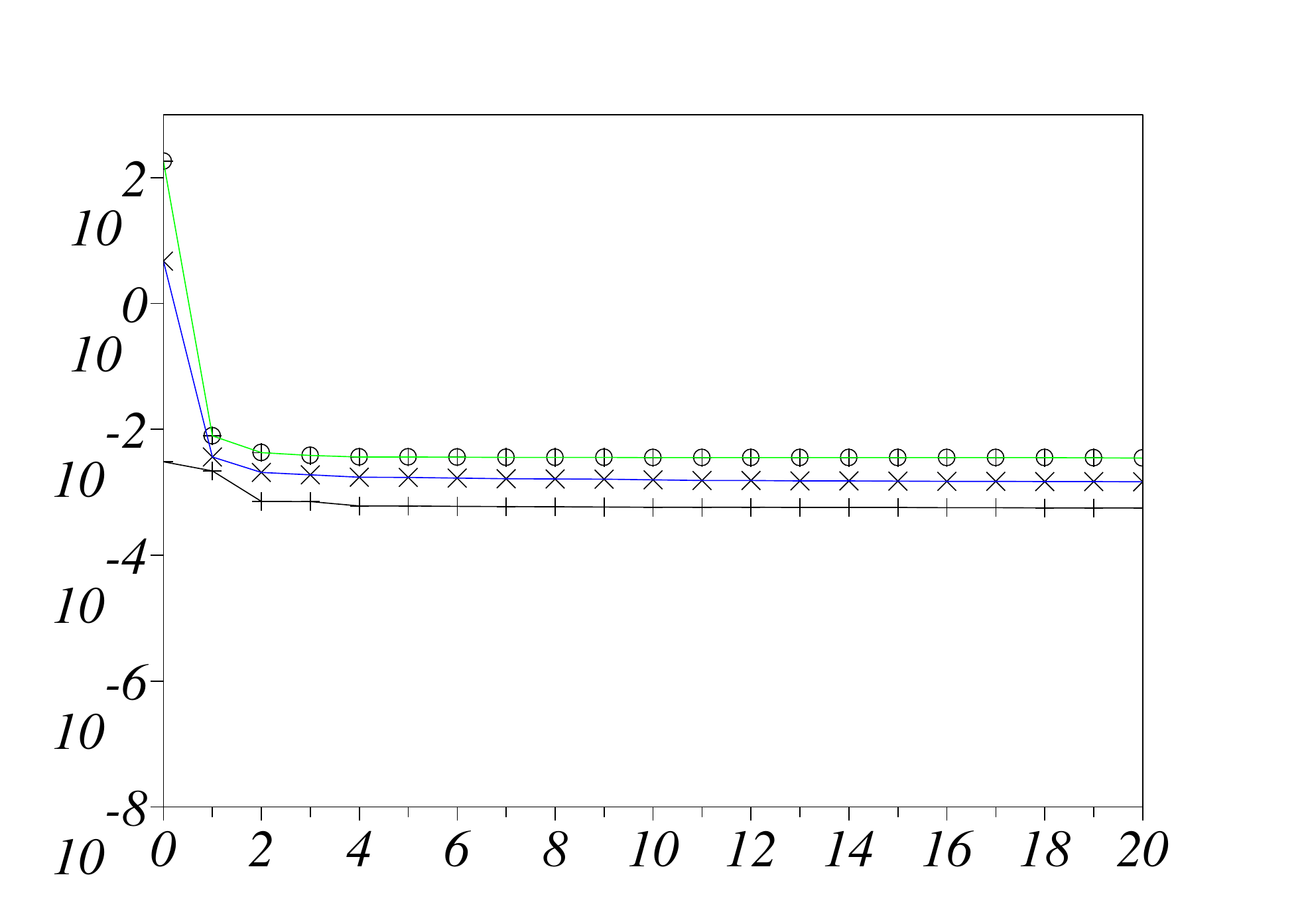}
\\
\includegraphics[trim = 8mm 0mm 20mm 10mm, clip, scale=.37]{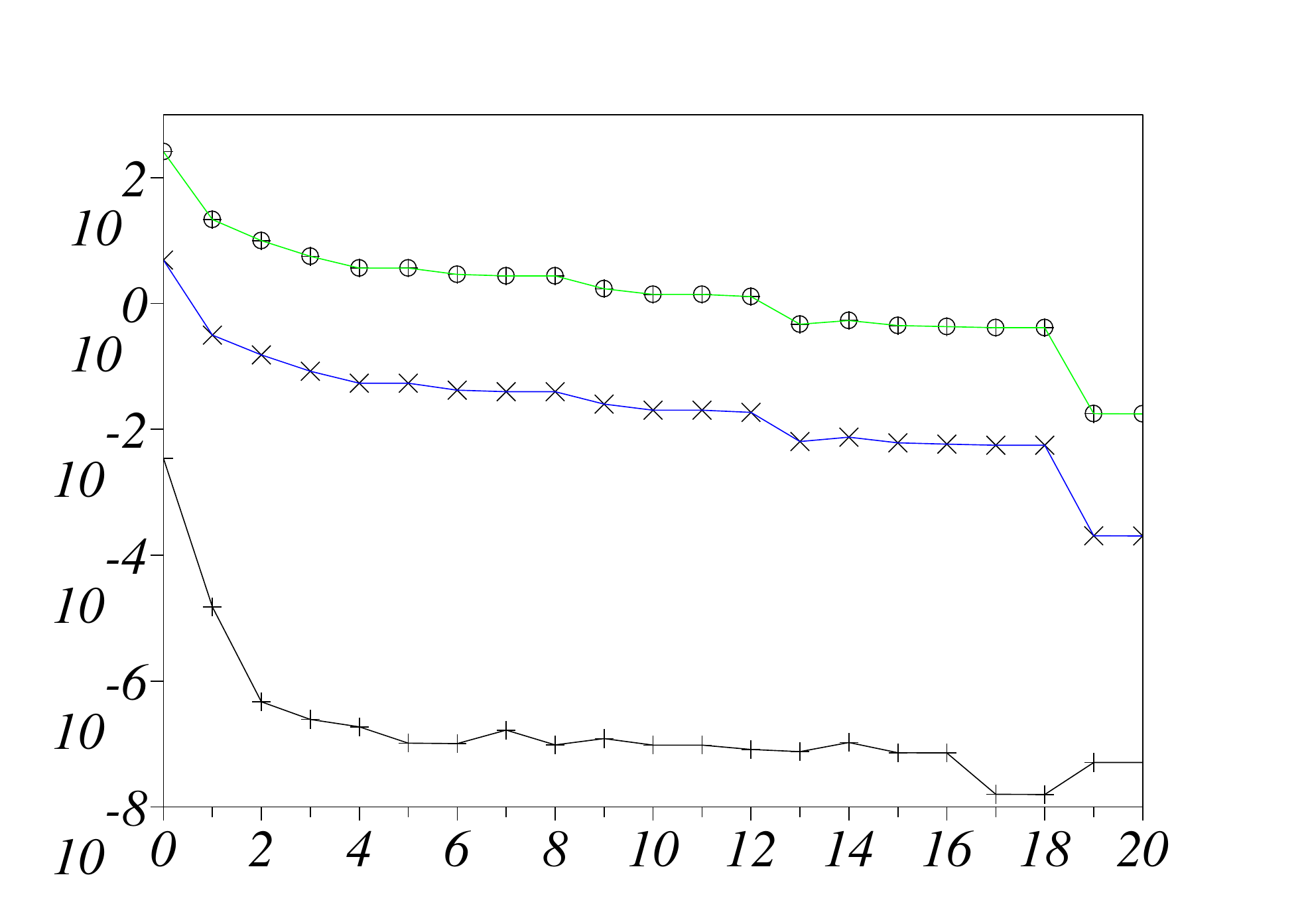}
\includegraphics[trim = 8mm 0mm 20mm 10mm, clip, scale=.37]{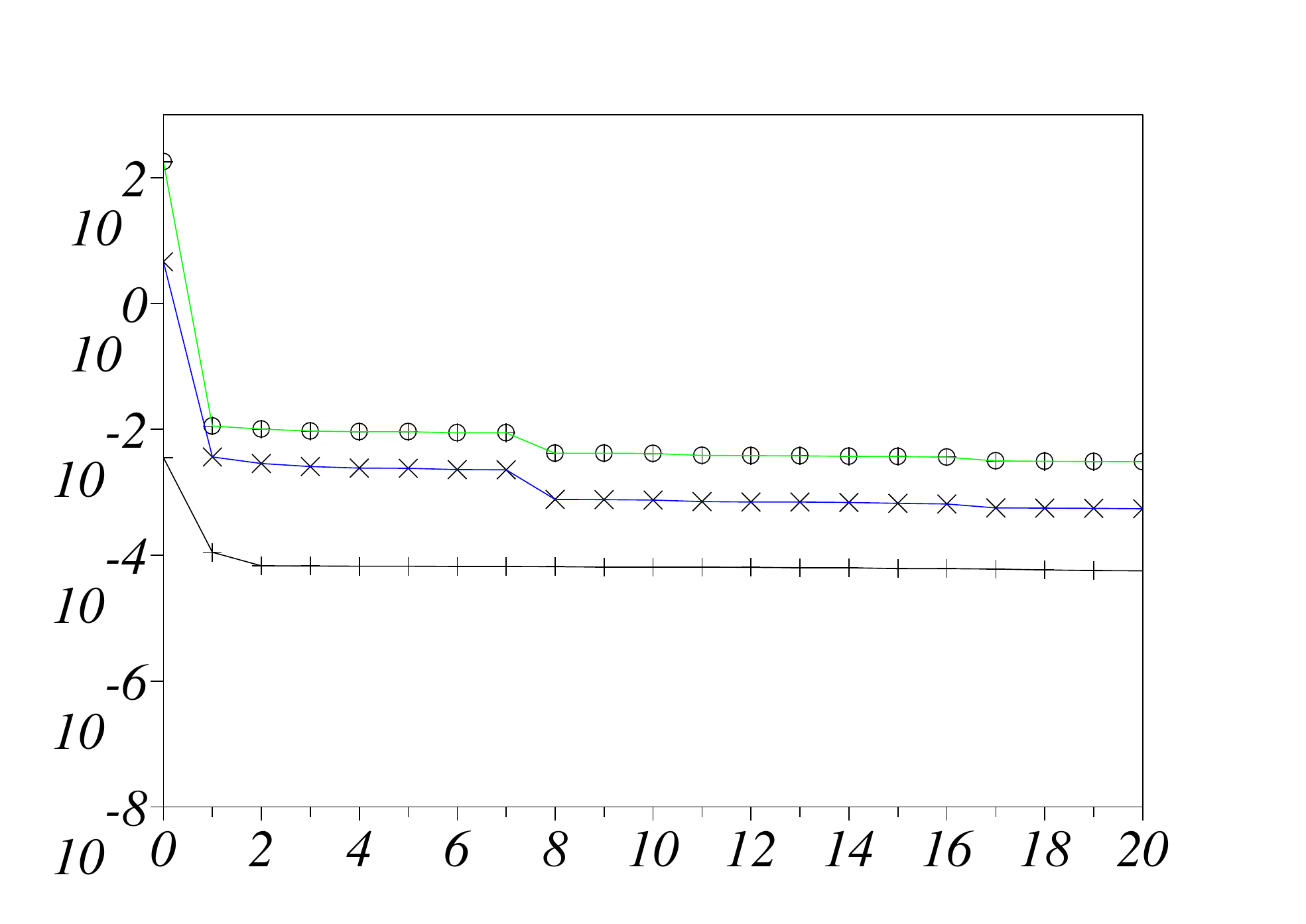}
\caption{Algorithm~1 (left) and 2 (right) for Black--Scholes model with local ``hyperbolic'' volatility: 
Minimum $+$, mean $\times$ and maximum $\circ$ of the 
(online) relative variance~\eqref{eq:relativerv} 
in a sample test $\Lambda_{\rm test wide}$ of parameters
with respect to the size $I$ of the reduced basis.
Greedy selection with absolute variance~\eqref{eq:absoluterv} (top) and 
relative variance~\eqref{eq:relativerv} (bottom).
}
\label{fig:bsdistrib_Erel_large}
\end{figure}

Now, in Fig.~\ref{fig:bsdistrib_large} and~\ref{fig:bsdistrib_Erel_large}, 
we show the online (absolute and relative) variance
for a new sample test of parameters $\Lambda_{\rm test wide}$ 
uniformly distributed in  
$\Lambda_{\rm wide}=[-.15,.25]\times\{b=c\in]0,2[\}\times\{1.\}\times\{1.1\}\times\{5\}\times\{.05\}$,
which is twice larger than 
$\Lambda=[-.05,.15]\times\{b=c\in[.5,1.5]\}\times\{1.\}\times\{1.1\}\times\{5\}\times\{.05\}$
where the training sample $\Lambda_{\rm trial}$ of the offline stage is nested~:
the quality of the variance reduction compared to that for
a narrower sample test $\Lambda_{\rm test}$ seems to decrease faster for
Algorithm~1 than for Algorithm~2.
So Algorithm~2 definitely seems more robust with respect to the variations in $\lambda$
than Algorithm~1.
This observation is even further increased 
if we use the relative variance~\eqref{eq:relativerv}
instead of the absolute variance~\eqref{eq:absoluterv}, 
as shown by the results in Fig.~\ref{fig:bsdistrib_large} and~\ref{fig:bsdistrib_Erel_large}.


\subsection{Vector processes with constant diffusion and parametrized drift} 
\label{sec:polymer}

\subsubsection{Molecular simulation of dumbbells in polymeric fluids}

In rheology of polymeric viscoelastic fluids,
the long polymer molecules responsible for the viscoelastic behaviour
can be modelled through kinetic theories of statistical physics as Rouse chains,
that is as chains of Brownian beads connected by springs.
We concentrate on the most simple of those models,
namely ``dumbbells'' (two beads connected by one spring)
diluted in a Newtonian fluid.

Kinetic models consist in adding to the usual velocity and pressure fields $(\uu{u},p)$
describing the (macroscopic) state of the Newtonian solvent,
a field of dumbbells represented by their end-to-end vector $\uu{X}_t(\uline{x})$ at time $t$
and position $\uline{x}$ in the fluid.
Vector stochastic processes $(\uu{X}_t(\uline{x}))$ encode the time evolution of
the orientation and the stretch of the dumbbells (the idealized configuration of a polymer molecule)
for each position $\uline{x}\in\mathcal{D}$
in a macroscopic domain $\mathcal{D}$ where the fluid flows.
To compute the flow of a viscoelastic fluid with such multiscale dumbbell models~\cite{le-bris-lelievre-07},
segregated algorithms are used that iteratively,
on successive time steps with duration $T$:
\begin{itemize}
 \item first evolve the velocity and pressure fields $(\uu{u},p)$
of the Newtonian solvent under 
a fixed extra (polymeric) stress tensor field $\uu{\tau}$
(typically following Navier-Stokes'equations), and
 \item then evolve the (probability distribution of the) polymer configurations vector field
 $(\uu{X}_t(\uline{x}))$ surrounded by the newly computed fixed 
velocity field $\uu{u}$.
\end{itemize}

The physics of kinetic models is based on a scale separation between 
the polymer molecules and the surrounding Newtonian fluid solvent.
On the one side, the polymer configurations are directly influenced by the (local)
velocity and pressure of the Newtonian solvent in which they are diluted.
Reciprocally, on the other side,
one needs to compute at every $\uline{x} \in \mathcal{D}$
the extra (polymeric) stress, given the Kramers formula:
$$ \uu{\tau}(T,\uline{x})=\E{ \uu{X}_T(\uline{x})\otimes\uu{F}(\uu{X}_T(\uline{x})) } \,,$$
after one evolution step $t\in[0,T]$ over which the polymer configurations have evolved
(remember that here $[0,T]$ should be understood as a timestep).
The vector valued process $\uu{X}_t(\uline{x})$ in $\R^d$ ($d=2$ or $3$) solves a Langevin equation
at every physical point $\uline{x} \in \mathcal{D}$ (Eulerian description):
\begin{equation*}
 d\uu{X}_t +\uu{u}\cdot\nabla_{\uline{x}} \uu{X}_t dt
= \left( (\nabla_{\uline{x}}\uu{u})\,\uu{X}_t-\uu{F}(\uu{X}_t) \right)\, dt + d\uu{B}_t \,.
\end{equation*}
This Langevin equation describes the evolution of polymers
at each $\uline{x} \in \mathcal{D}$,
under an advection $\uu{u}\cdot\nabla_{\uline{x}} \uu{X}_t$, 
a hydrodynamic force $(\nabla_{\uline{x}}\uu{u})\uu{X}_t$,
Brownian collisions $(\uu{B}_t)$ with the solvent molecules, and 
an entropic force $\uu{F}(\uu{X}_t)$ specific to the polymer molecules.
Typically, this entropic force reads either $\uu{F}(\uu{X})=\uu{X}$ (for Hookean dumbbells), 
or $\uu{F}(\uu{X})=\frac{\uu{X}}{1 - |\uu{X}|^2/b}$
(for Finitely-Extensible Nonlinear Elastic or FENE dumbells, to model the finite extensibility of polymers: $|\uu{X}|<\sqrt{b}$).

In the following, we do not consider the advection term $\uu{u}\cdot\nabla_{\uline{x}} \uu{X}_t$
(which can be handled through integration of the characteristics in a semi-Lagrangian framework, for instance),
and we concentrate on solving the parametrized SDE:
\begin{equation}   \label{eq:langevin}
 d\uu{X}_t = \left( \uuline{\lambda} \, \uu{X}_t - \uu{F}(\uu{X}_t) \right)\, dt + d\uu{B}_t \,,
\end{equation}
on a time slab $[0,T]$, with a fixed matrix $\uuline{\lambda}(\uline{x})=\nabla_x\uu{u}(\uline{x})$.
We also assume, as usual for viscoelastic fluids, that
the velocity field is incompressible (that is $\tr(\uuline{\lambda})=0$), 
hence the parameter $\uuline{\lambda}$ is only $(d^2-1)$-dimensional.

This is a typical many-query context where 
the Langevin equation~\eqref{eq:langevin} has to be computed many times at each (discretized) position
$\uline{x} \in \mathcal{D}$, for each value of the $d\times d$-dimensional parameter $\uuline{\lambda}$ (since $\nabla_{\uline{x}}\uu{u}(\uline{x})$ depends on the position $\uline{x}$).
Furthermore, the computation of the time-evolution of the flow
defines a {\it very demanding} many-query context where the latter has to be done iteratively
over numerous time steps of duration $T$ 
between which the tensor field $\uuline{\lambda}(\uline{x})$ is evolved
through a macroscopic equation for the velocity field $\uu{u}$.

\begin{rem}[Initial Condition of the SDE as additional parameter]
\label{rem:initial_condition}
Let $T_0=0$ and $T_{n+1} = (n+1) T$.
Segregated numerical schemes for kinetic models of polymeric fluids
as described above simulate~\eqref{eq:langevin}
on successive time slabs $[T_n,T_{n+1}]$, for $n\in\mathbb{N}$.
More precisely, on each time slab $[T_n,T_{n+1}]$, 
one has to compute 
\begin{align}\label{condexp}
\nonumber
\uu{\tau}(T_{n+1})
& = \E{ \uu{X}_{T_{n+1}}\otimes\uu{F}(\uu{X}_{T_{n+1}}) }
\\
& =\E{ \E{ \uu{X}_{T_{n+1}}\otimes\uu{F}(\uu{X}_{T_{n+1}})| \uu{X}_{T_{n}} } }
\end{align}
at a fixed position $\uline{x}\in\mathcal{D}$.
In practice,~\eqref{condexp} can be approximated through
\begin{equation}\label{eq:condexpapprox}
\uu{\tau}(T_{n+1}) 
\simeq \frac1R \sum_{r=1}^R \frac1M \sum_{m=1}^M 
\uu{X}_{T_{n+1}}^{r,m}\otimes\uu{F}(\uu{X}_{T_{n+1}}^{r,m}) \,,
\end{equation}
after simulating $M R$ processes $(\uu{X}_{t}^{r,m})_{t\in[T_n,T_{n+1}]}$
driven by $M R$ independent Brownian motions
for a given set of $R$ {\it different} initial conditions, typically:
$$ \uu{X}_{T_n^+}^{r,m} = \uu{X}_{T_n^-}^{r,1}\,,\ r=1,\ldots,R\,,\ m=1,\ldots,M\,, $$
or any $\uu{X}_{T_n^-}^{r,m_0}$ (with $1\le m_0\le M$) given by
the computation at final time of the previous time slab $[T_{n-1},T_n]$.
In view of~\eqref{eq:condexpapprox}, 
for a fixed $r$, the computation of $\frac1M \sum_{m=1}^M 
\uu{X}_{T_{n+1}}^{r,m}\otimes\uu{F}(\uu{X}_{T_{n+1}}^{r,m})$
using the algorithms presented above requires a modification of the methods to the case when the initial condition of the SDE assumes many values.

To adapt Algorithm~1 to the context of SDEs with many different initial conditions,
one should consider reduced bases 
for control variates 
which depend on the 
joint-parameter $({\lambda},x)$ where $x$ is the initial condition of the SDE.
And variations on the joint-parameter $({\lambda},x)$
can be simply recast into the framework of
SDEs with fixed initial condition 
used for presentation of Algorithm~1
after the change of variable $\hat{X}^{\lambda,x}_t={X}^{\lambda}_t-x$,
that is using the family of SDEs
with fixed initial condition ${\hat{X}^{\lambda,x}}_0=0$:
\begin{equation}\label{eq:langevinbis}
d{\hat{X}}^{\lambda,x}_t 
= \hat{b}^{\lambda,x}(t,\hat{X}^{\lambda,x}_t) dt 
+ \hat{\sigma}^{\lambda,x}(t,\hat{X}^{\lambda,x}_t) d{B}_t \,,
\end{equation}
where $\hat{b}^{\lambda,x}(t,X)=b^{\lambda}(t,X+x)$,
$\hat{\sigma}^{\lambda,x}(t,X)=\sigma^{\lambda}(t,X+x)$,
for all $t$, $X$ and $x$.
Then, with $\hat{g}^{\lambda,x}(X)=g^{\lambda}(X+x)$ and
$\hat{f}^{\lambda,x}(t,X)=f^{\lambda}(t,X+x)$, 
the output is the expectation of
$$ 
\hat{Z}^{\lambda,x}=\hat{g}^\lambda(\hat{X}_T^{\lambda,x})
-\int_0^T \hat{f}^\lambda(s,\hat{X}_s^{\lambda,x})\, ds \,.
$$
And the corresponding ``ideal'' control variate reads 
$ 
\hat{Y}^{\lambda,x}
= \hat{Z}^{\lambda,x}
-\E{\hat{Z}^{\lambda,x}} \,.
$ 

In Algorithm~2,
note that $u^{\lambda}$ solution to~\eqref{eq:PDE}
does not depend on the initial condition used for the SDE.
So, once parameters ${\lambda_i}$ ($i=1,\ldots,I$) have been selected offline,
Algorithm~2 applies similarly for SDEs with one fixed, or many different, initial conditions.
Though, the offline selection of parameters ${\lambda_i}$
using SDEs with many different initial conditions
should consider a larger trial sample than for one fixed initial condition.
Indeed, the selection criterium in the greedy algorithm does depend on the initial condition of the SDE.
So, defining a trial sample of initial conditions $\Lambda_{\rm IC}$,
the following selection should be performed at step $i$ in Fig.~\ref{fig:greedy2}:
$$  
\textnormal{Select }
\lambda_{i+1} \in 
\underset{\lambda \in\Lambda_{\rm trial}\backslash\{\lambda_j,j=1,\ldots, i\}}{\argmax}
\underset{x \in\Lambda_{\rm IC}}{\max}
{\rm Var_{M_{small}}}({Z}^{\lambda,x} -{\tilde{Y}^{\lambda,x}_i}) 
\,,
$$
where ${Z}^{\lambda,x}$ and ${\tilde{Y}^{\lambda,x}_i}$,
defined like $Z^\lambda$ and $\tilde{Y}^\lambda_i$,
depend on $x$ because the stochastic process $(X^\lambda_t)$ depends on $X^\lambda_0=x$.

It might be useful to build different reduced bases,
one for each cell of a partition of the set of the initial condition.
In summary, both algorithms can be extended to SDEs with variable initial condition,
at the price of increasing the dimension of the parameter
(see also Remark~\ref{rem:highdim}).
\end{rem}

\begin{rem}[Multi-dimensional output]\label{rem:output}
Clearly, the full output $\uu{\tau}$ in the problem described above is three-dimensional
(it is a symmetric matrix).
So our reduced-basis approach such as presented so far
would need three different reduced bases,
one for each scalar output.
Though, one could alternatively consider
the construction of only one reduced basis for the three outputs,
which may be advantageous, 
see~\cite{boyaval-08} for one example of such a construction.
\end{rem}

Note that it is difficult to compute accurate approximations
of the solution to the backward Kolmogorov equation~\eqref{eq:PDE}
in the FENE case, because of the nonlinear explosive term.
It is tractable in some situations,
see~\cite{lozinski-chauviere-03,chauviere-lozinski-03} for instance,
though at the price of computational difficulties we did not want to 
deal with in this first work on our new variance reduction approach.
On the contrary, the backward Kolmogorov equation~\eqref{eq:PDE}
can be solved exactly in the case of Hookean dumbells.
Hence we have approximated here $u^\lambda$ in Algorithm~2
by the numerical solution $\tilde{u}^\lambda$ to
the backward Kolmogorov equation~\eqref{eq:PDE} for Hookean dumbells,
whatever the type of dumbbells used for the molecular simulation
(Hookean or FENE).

We would like to mention the recent work~\cite{knezevic-patera-09}
where the classical reduced-basis method for parameterized PDEs has been used in the FENE case (solving the FENE Fokker-Planck by dedicated deterministic methods).
Our approach is different since we consider 
a stochastic discretization.

 

\subsubsection{Numerical results}


The SDE~\eqref{eq:pb} for FENE dumbbells (when $d=2$)
is discretized with the Euler-Maruyama scheme
using $N=100$ iterations with a constant time step of $\Delta t = 10^{-2}$
starting from a (deterministic) initial condition $\uu{X}_0=(1,1)$,
with reflecting boundary conditions at the boundary of the ball with radius $\sqrt{b}$.

The number of realizations used for the Monte-Carlo evaluations,
and the sizes of the (offline) trial sample $\Lambda_{\rm trial}$ 
and (online) test sample $\Lambda_{\rm test}$
for the three-dimensional matrix parameter $\uuline{\lambda}$
with entries $(\lambda_{11}=-\lambda_{22},\lambda_{12},\lambda_{21})$,
are kept similar to the previous Section~\ref{sec:bs}.
Samples $\Lambda_{\rm trial}$ and $\Lambda_{\rm test}$ for the parameter $\uuline{\lambda}$ are uniformly distributed in a cubic range $\Lambda= [-1,1]^3$.
We will also make use of an enlarged (online) test sample $\Lambda_{\rm test wide}$,
uniformly distributed in the range $[-2,2]^3$.

\begin{figure}[htbp]
\includegraphics[trim = 8mm 0mm 20mm 10mm, clip, scale=.37]{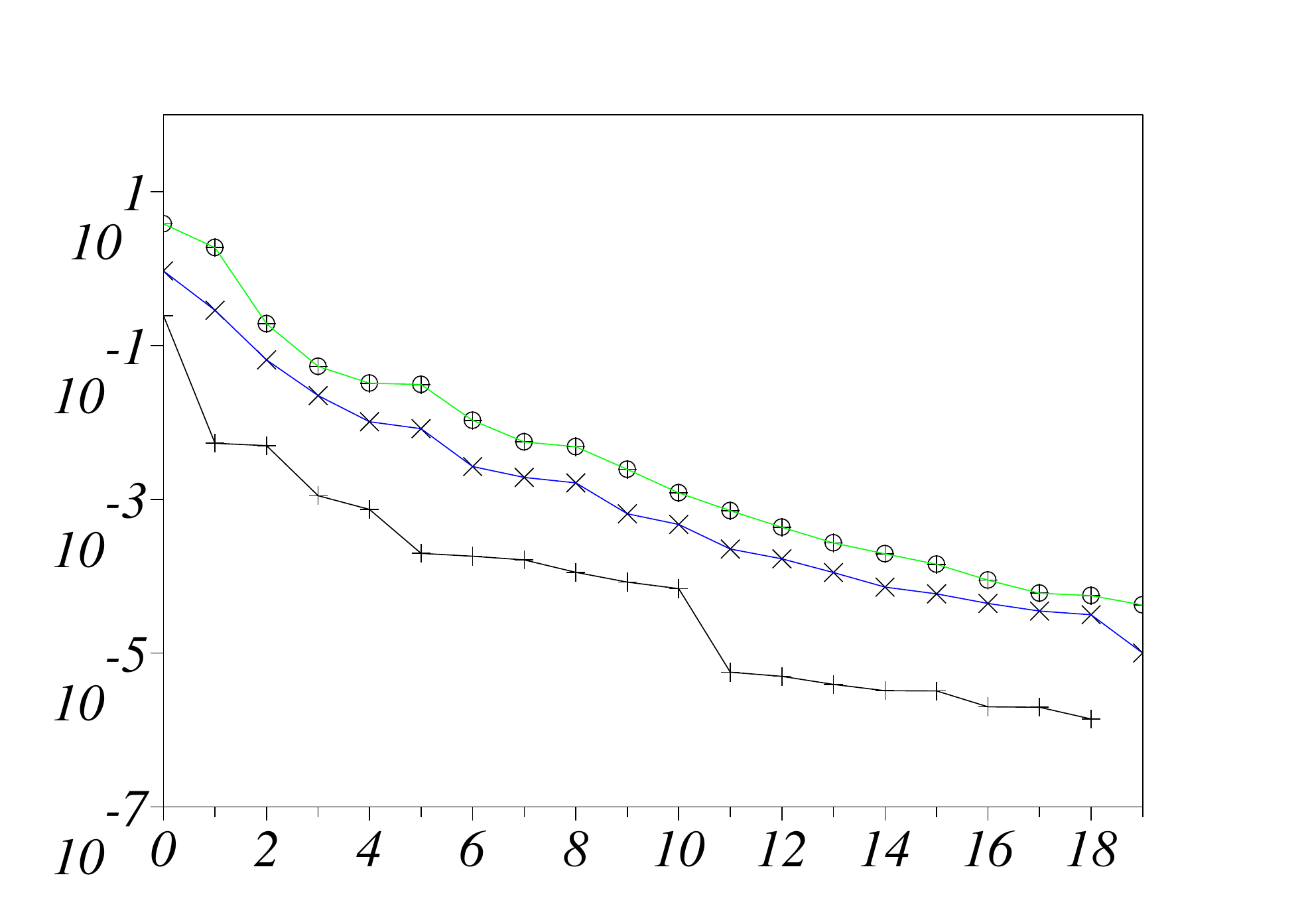}
\includegraphics[trim = 8mm 0mm 20mm 10mm, clip, scale=.37]{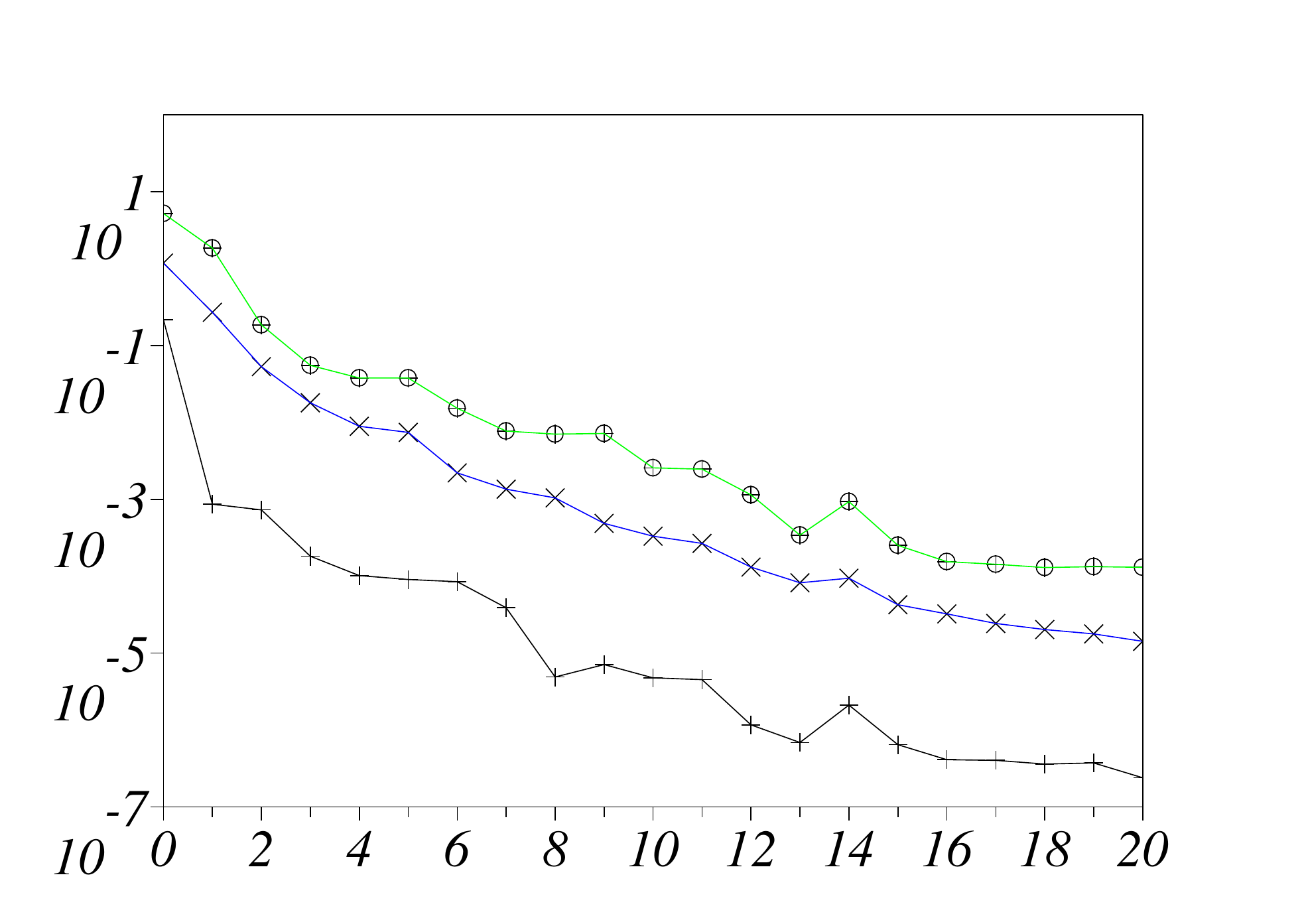}
\caption{Algorithm~1 for FENE model with $b=9$: 
Minimum $+$, mean $\times$ and maximum $\circ$ of the absolute variance~\eqref{eq:absoluterv} 
in samples of parameters
(left: offline sample $\Lambda_{\rm trial}\setminus\{\lambda_i,i=1,\ldots,I\}$; right: online sample $\Lambda_{\rm test}$) 
with respect to the size $I$ of the reduced basis.
\label{fig:FENEdistrib1}
}
\includegraphics[trim = 8mm 0mm 20mm 10mm, clip, scale=.37]{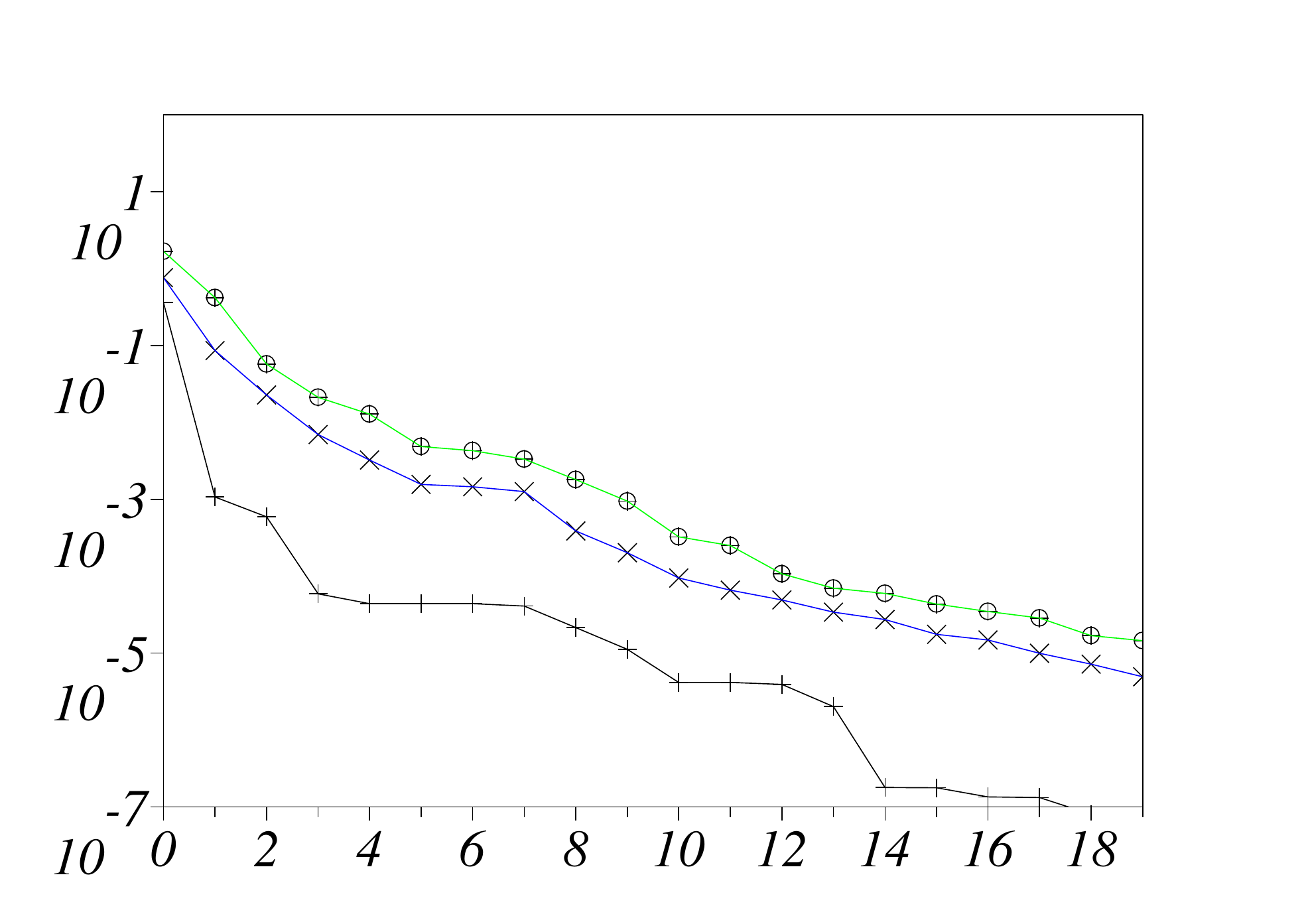}
\includegraphics[trim = 8mm 0mm 20mm 10mm, clip, scale=.37]{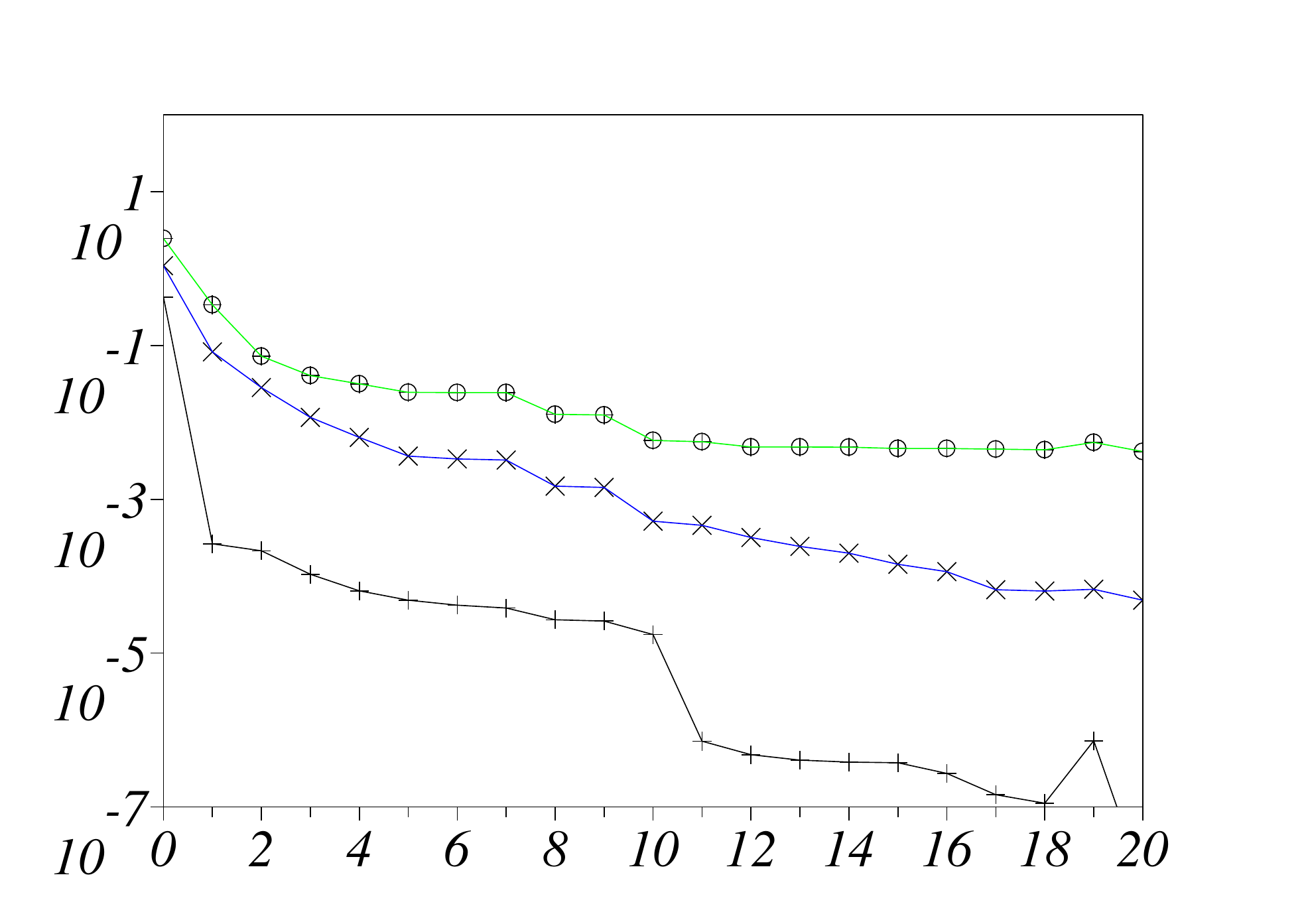}
\caption{Algorithm~1 for FENE model with $b=4$: 
Minimum $+$, mean $\times$ and maximum $\circ$ of the absolute variance~\eqref{eq:absoluterv} 
in samples of parameters
(left: offline sample $\Lambda_{\rm trial}\setminus\{\lambda_i,i=1,\ldots,I\}$; right: online sample $\Lambda_{\rm test}$) 
with respect to the size $I$ of the reduced basis.
\label{fig:FENEdistrib3}
}
\includegraphics[trim = 8mm 0mm 20mm 10mm, clip, scale=.37]{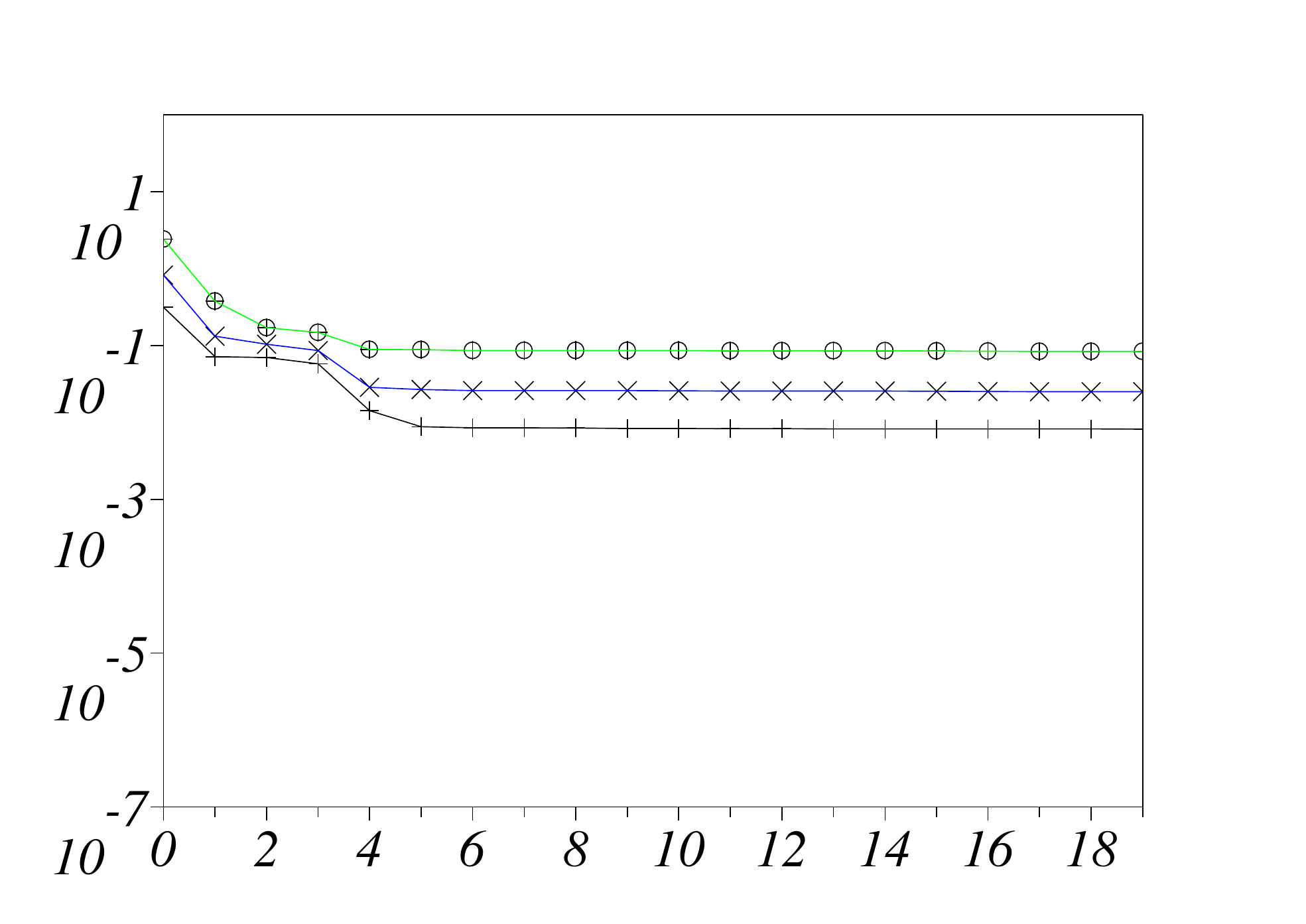}
\includegraphics[trim = 8mm 0mm 20mm 10mm, clip, scale=.37]{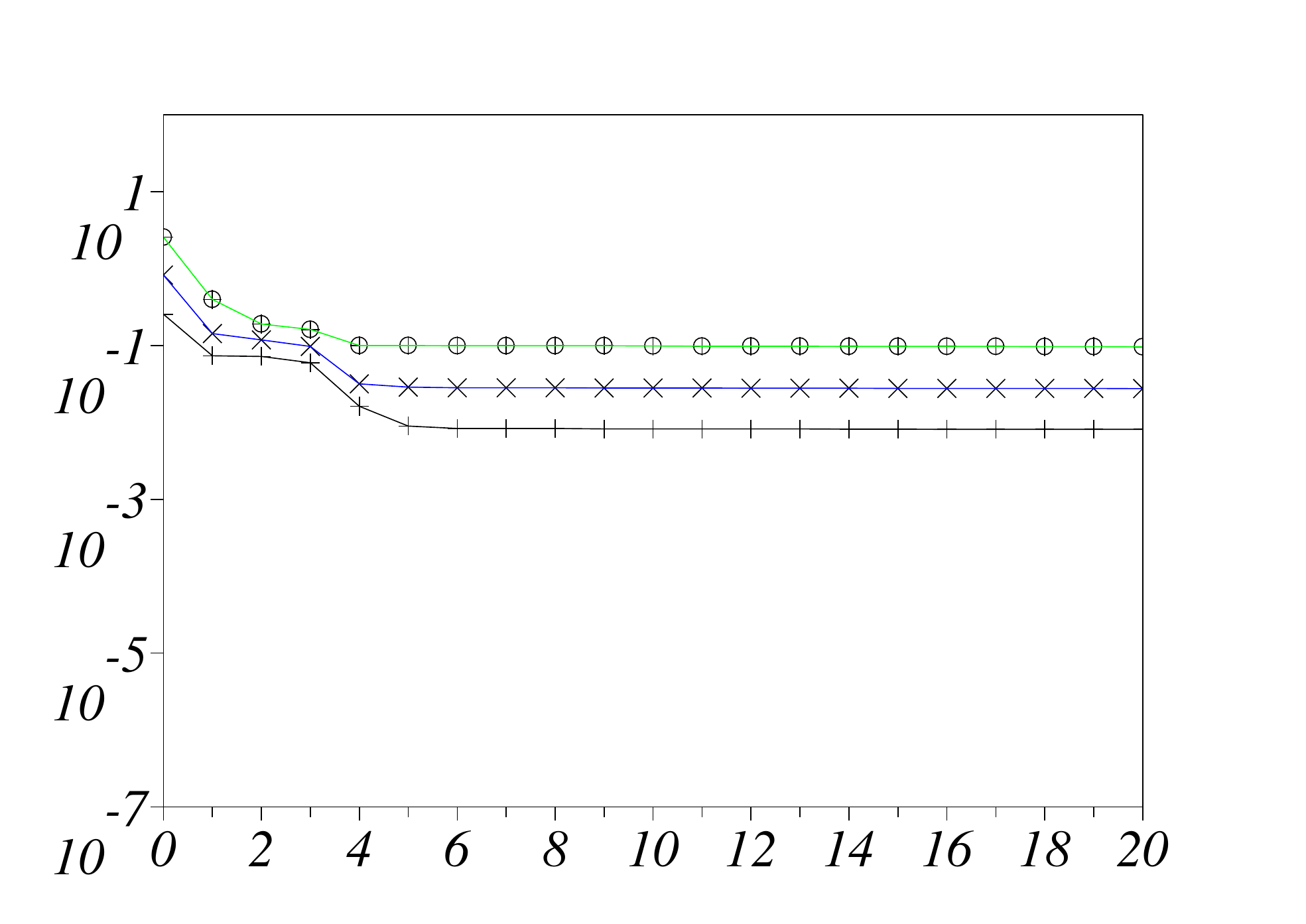}
\caption{Algorithm~2 for FENE model with $b=9$: 
Minimum $+$, mean $\times$ and maximum $\circ$ of the absolute variance~\eqref{eq:absoluterv} 
in samples of parameters
(left: offline sample $\Lambda_{\rm trial}\setminus\{\lambda_i,i=1,\ldots,I\}$; right: online sample $\Lambda_{\rm test}$) 
with respect to the size $I$ of the reduced basis.
\label{fig:FENEdistrib5}
}
\end{figure}

When $b=9$, the variance reduction online with Algorithm~1 is again very interesting,
of about $4$ orders of magnitude with $I=20$ basis functions,
whatever the criterium used for the selection
(we only show the absolute variance, in Fig.~\ref{fig:FENEdistrib1}).
But when $b=4$, the reflecting boundary conditions are more often active,
and the maximum online variance reduction slightly degrades 
(see Fig.~\ref{fig:FENEdistrib3}).

\begin{figure}[htbp]
\includegraphics[trim = 8mm 0mm 20mm 10mm, clip, scale=.37]{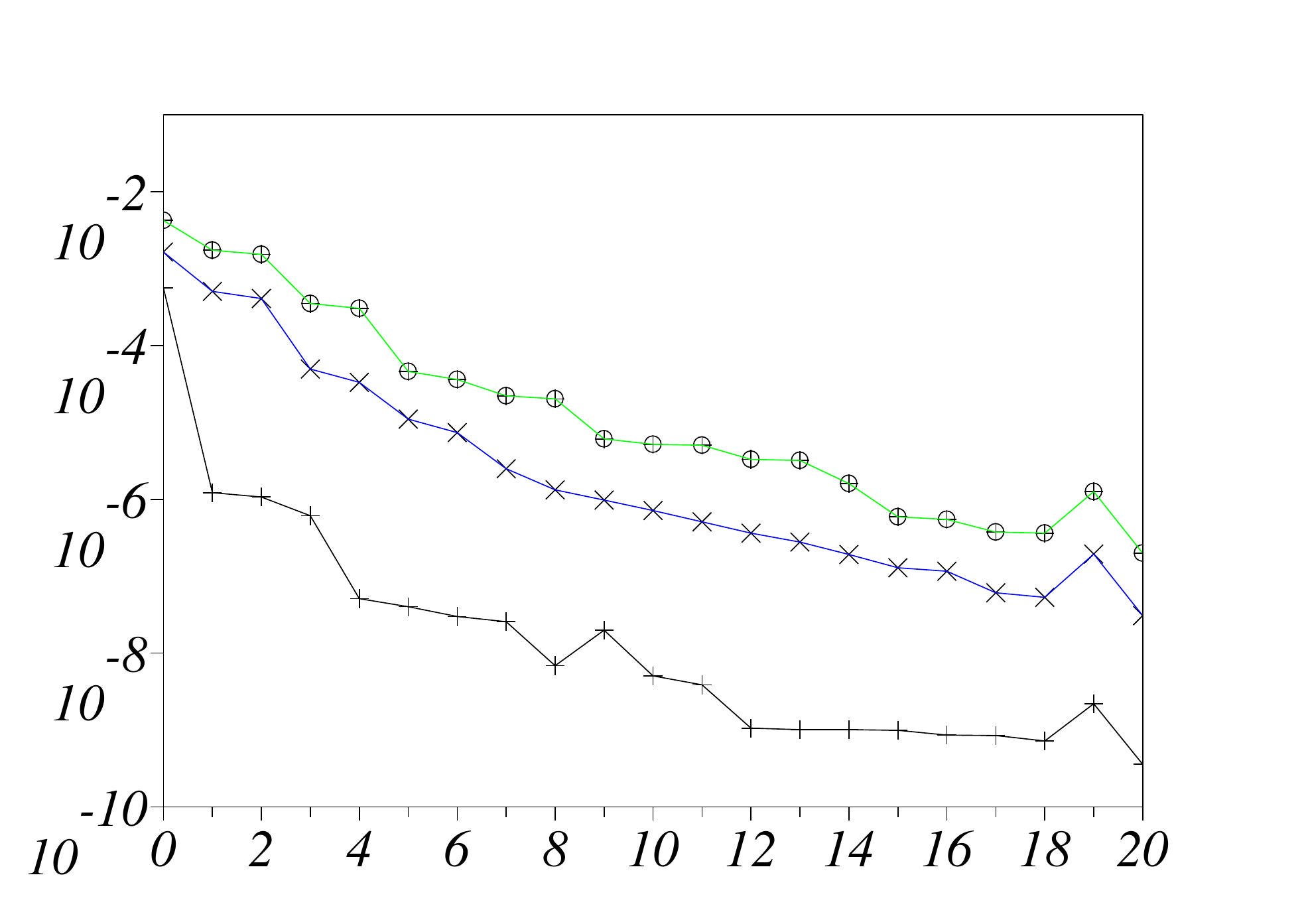}
\includegraphics[trim = 8mm 0mm 20mm 10mm, clip, scale=.37]{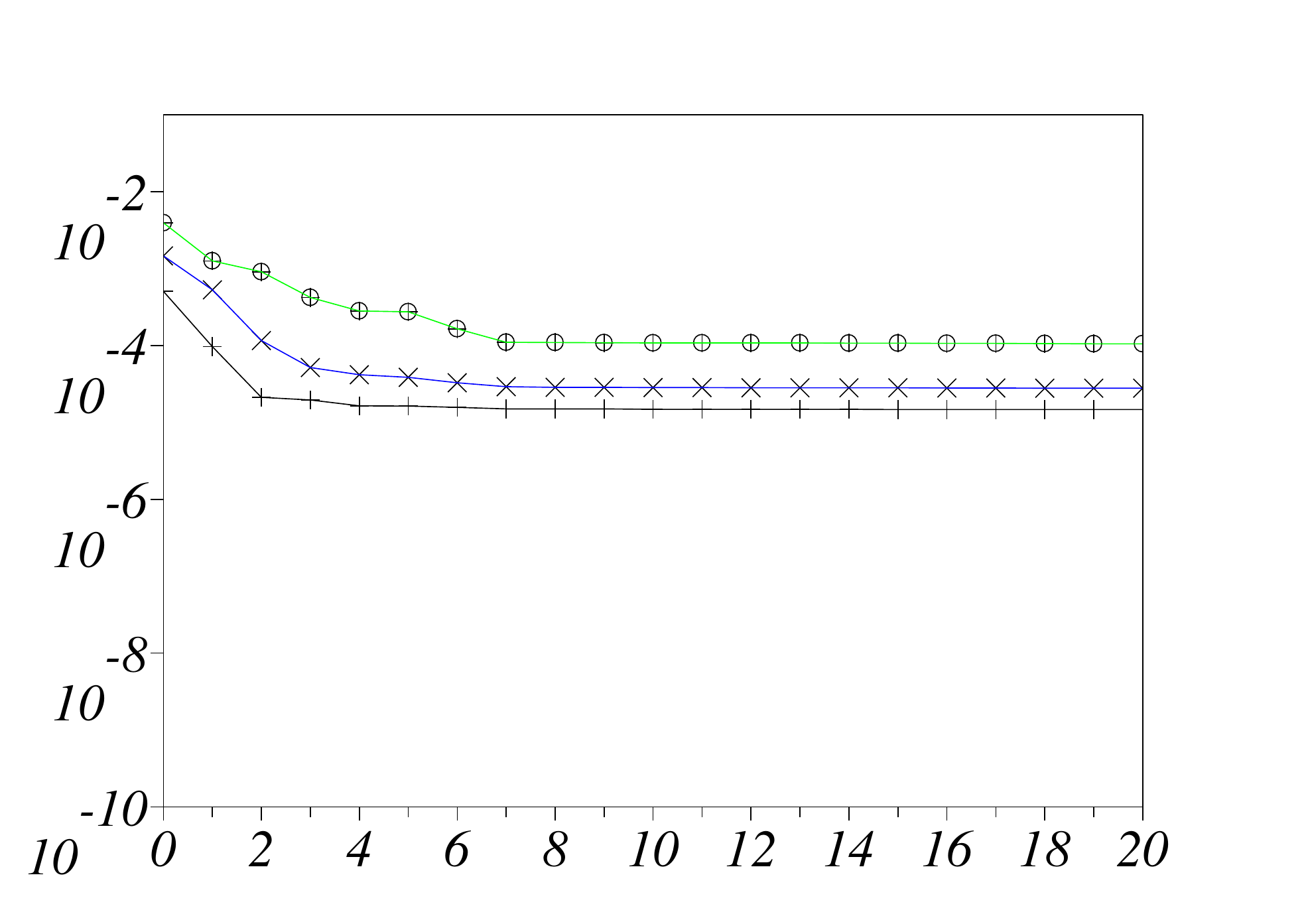}

\includegraphics[trim = 8mm 0mm 20mm 10mm, clip, scale=.37]{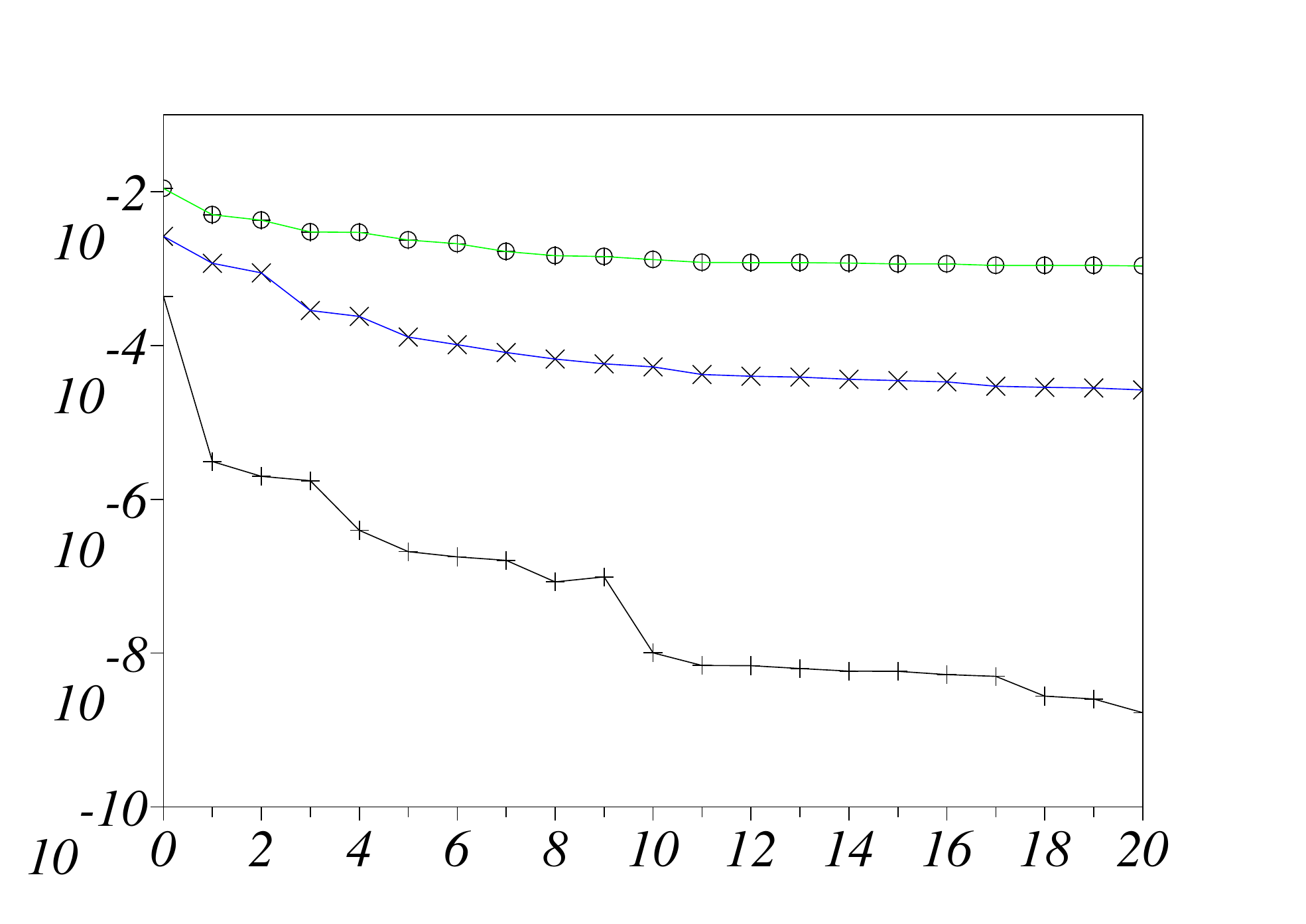}
\includegraphics[trim = 8mm 0mm 20mm 10mm, clip, scale=.37]{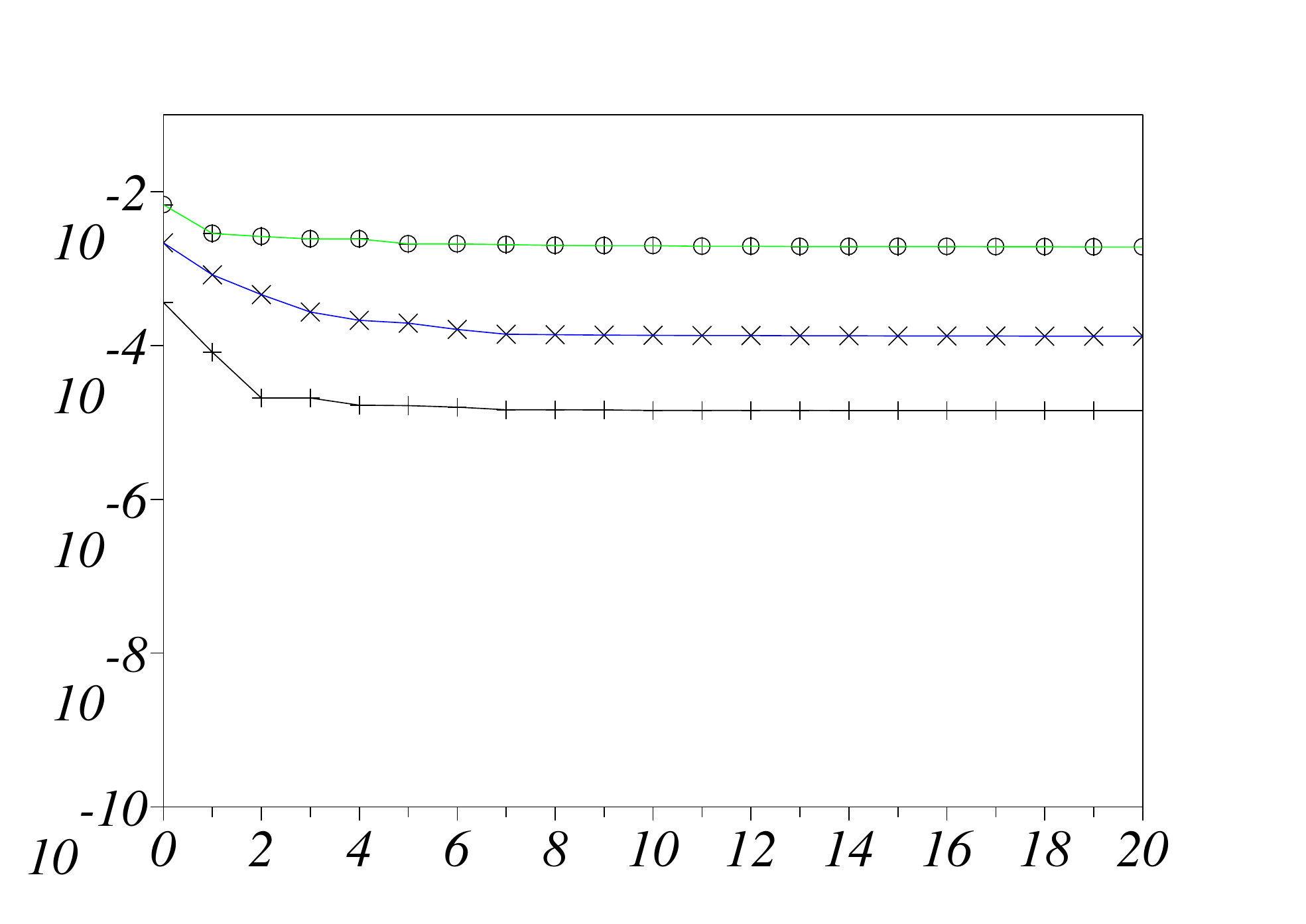}
\caption{Algorithm~1 (left) and 2 (right) for FENE model with $b=16$:
Minimum $+$, mean $\times$ and maximum $\circ$ of the relative variance~\eqref{eq:relativerv}
in online test for samples $\Lambda_{\rm test}$ (top)
and $\Lambda_{\rm test wide}$ (bottom) of parameters,
with respect to the size $I$ of the reduced basis.
\label{fig:FENEdistrib6}
}
\end{figure}

We first tested our variance reduction with Algorithm~2 for Hookean dumbells
and it appeared to work well;
but such a model is considered too simple generally.
Then using the solution to the Kolmogorov backward equation for Hookean dumbells
as $\tilde u^\lambda$
in Algorithm~2 for FENE dumbbells still yields good variance reduction 
while the boundary is not touched 
(see Fig.~\ref{fig:FENEdistrib5}) ;
when $b=4$ and many reflections at the boundary occur,
the variance is hardly reduced.
Again Algorithm~2 seems to be slightly more robust than Algorithm~1
in terms of extrapolation,
that is when the (online) test sample
is ``enlarged''
(see Fig.\ref{fig:FENEdistrib6} with $b=16$ and
a sample test (online) $\Lambda_{\rm test wide}$).


\section{Conclusion and perspectives}
\label{sec:ccl}

We have demonstrated the feasibility of a reduced-basis approach 
to compute control variates 
for the expectation of functionals of a parameterized It\^o stochastic process.
We have also tested the efficiency of such an approach with two possible algorithms,
in two simple test cases where either the drift or the diffusion of
scalar ($d=1$), and vector ($d=2$), It\^o processes are parametrized, 
using $2$-~or $3$-dimensional parameters.

Algorithm~2 is less generic than Algorithm~1 ; 
it is basically restricted to low-dimensional stochastic processes $(X_t)$ since:
	\begin{itemize}
	\item one needs to solve (possibly high-dimensional) PDEs (offline), and
	\item discrete approximations of the PDEs solutions on a grid have to be kept 
	in memory (which is possibly a huge amount of data).
	\end{itemize}
Yet, Algorithm~2 seems more robust to variations in the parameter.

From a theoretical viewpoint,
it remains to better understand the convergence of 
reduced-basis approximations for parametrized control variates
depending on the parametrization
(and on the dimension of the parameter in particular),
on the reduced-basis construction
(following a greedy procedure)
and on an adequate discretization choice 
(including the computation of approximate control variates
and the choice of a trial sample $\Lambda_{\rm trial}$).

\medskip
{\bf Acknowledgement.}
We acknowledge financial support from the France Israel Teamwork in Sciences.
We thank Raz Kupferman,
Claude Le Bris, Yvon Maday and Anthony T. Patera for fruitful discussions.
We are grateful to the referees for constructive remarks.
\medskip

\appendix


\section{\emph{\ Algorithm~2 in a higher-dimensional setting ($d\ge4$)}}
\label{app:pde_discretization}

\setcounter{section}{6}

The solution $u^\lambda(t,y)$ to~\eqref{eq:PDE} 
can be computed at any $(t,y)\in [0,T] \times \R^d$
by the martingale representation theorem~\cite{karatzas-shreve-91}:
\begin{equation}\label{martingale_representation}
g^{\lambda}(X_T^{\lambda})-\int_t^T f^{\lambda}(s,X_s^{\lambda}) ds
= u^\lambda(t,X^{\lambda}_t) + 
\int_t^T \nabla u^\lambda(s,X_s^{\lambda}) \cdot \sigma^{\lambda}(s,X_s^{\lambda}) dB_s,
\end{equation}
obtained by an It\^{o} formula similar to~\eqref{integral_representation}.
This gives the following Feynman-Kac formula for $u^\lambda(t,x)$,
which can consequently be computed at any $(t,y)\in [0,T] \times \R^d$ 
through Monte-Carlo evaluations:
\begin{equation}\label{eq:representation}
u^{\lambda}(t,y)
=\E{ g^{\lambda}(X_T^{\lambda,t,y})-\int_{t}^T f^{\lambda}(s,X_{s}^{\lambda,t,y}) ds },
\end{equation}
where~$(X_t^{\lambda,t_0,y})_{t_0 \le t \le T}$ is the solution to~\eqref{eq:pb} with initial condition $X_{t_0}^{\lambda,t_0,y}=y$. Differentiating~\eqref{eq:representation} (provided $f^{\lambda}$
and $g^{\lambda}$ are differentiable),
we even directly get a Feynman-Kac formula for $\nabla u^\lambda(t,y)$:
\begin{equation}\label{eq:D_prob_rep}
\nabla u^{\lambda}(t,y)
= \E{ \Phi_{T}^{\lambda,t,y} \cdot \nabla g^{\lambda}(X_T^{\lambda,t,y})
- \int_{t}^T \Phi_{s}^{\lambda,t,y} \cdot \nabla f^{\lambda}(s,X_{s}^{\lambda,t,y}) ds }
\end{equation}
where the stochastic processes 
$\left( \Phi_{s}^{\lambda,t,y} , s \in [t,T] \right)$ in $\R^{d\times d}$ satisfy 
the first-order variation of the SDE~\eqref{eq:pb} with respect to the initial condition,
that is $ \Phi_{s}^{\lambda,t,y} = \nabla_y X_s^{\lambda,t,y}$ for any $s \in [t,T]$:
\begin{equation}\label{eq:D_X_x_0}
\Phi_{s}^{\lambda,t,y} = \I_d + 
\int_{t}^s \Phi_{s'}^{\lambda,t,y} \cdot \nabla b^\lambda(s',X_{s'}^{\lambda,t,y}) ds'
+ \int_{t}^s \Phi_{s'}^{\lambda,t,y} \cdot \nabla \sigma^\lambda(s',X_{s'}^{\lambda,t,y}) dB_{s'},
\end{equation}
where $\I_d$ denotes the ${d\times d}$ identity matrix
(see~\cite{newton-94} for a more general and rigorous presentation of this Feynman-Kac formula 
in terms of the Malliavin gradient).
The stochastic integral~\eqref{eq:ideal_control} can then be computed for each realization of $(B_t)$,
after discretizing $\left( \Phi_{s}^{\lambda,t,y} , s \in [t,T] \right)$.

Discrete approximations of the Feynman-Kac formula~\eqref{eq:D_prob_rep}
have already been used succesfully in the context of computing control variates
for the reduction of variance, in~\cite{newton-94} for instance.
Note that 
this numerical strategy to compute $\nabla u^\lambda$
from a Feynman-Kac formula requires a lot of computations.
Yet, most often, the computation time of the functions $(t,y)\to\nabla u^\lambda(t,y)$ 
would not be a major issue in a reduced-basis approach, 
since this would be done \textit{offline} 
(that is, in a pre-computation step, once for all)
for only a few selected values of the parameter $\lambda$.
What is nevertheless necessary for the reduced-basis approach to work is the possibility
to store the big amount of data corresponding to a discretization of $\nabla u^\lambda(t,y)$ 
on a grid for the variable $(t,y)\in[0,T]\times\R^d$,
(the parameter $\lambda$ then assuming only a few values in $\Lambda$ 
--- of order $10$ in our numerical experiments ---),
and to have rapid access to those data in the \textit{online} stage
(where control variates are computed for any $\lambda\in\Lambda$ using those precomputed data).

\setcounter{section}{1}
 
\section{\emph{\ Proof of Proposition~\ref{prop:apriori_rb}}}
\label{app:apriori_rb}

\setcounter{section}{7}

First note that, since $\E{Y^\lambda}=0$, then $\Var{Z^\lambda}=\Var{Y^\lambda}$.

So, for all $\lambda\in\Lambda$ and 
for every linear combination of $Y^{\lambda^N_n}$, $n=1,\ldots,N$~: 
$$ Y_N = \sum_{n=1}^N a_n(\lambda) \sum_{j=1}^J g_j(\lambda_n^N)Y_j $$
(with any choice $a_n(\lambda)\in\R$, $\lambda_n^N\in\Lambda$, $n=1,\ldots,N$), there holds (recall that the $Y_j$, $j=1,\ldots,J$, are uncorrelated)~:
\begin{align}\nonumber
\Var{Z^\lambda-Y_N} 
& = \Var{Y^\lambda-Y_N} 
 = \int_\Omega \left|
\sum_{j=1}^J 
\left( g_j(\lambda)-\sum_{n=1}^N a_n(\lambda) g_j(\lambda_n^N)\right)  \, Y_j \right|^2 d\P
\\ \label{eq:rb1}
& \le \left(\sum_{j=1}^J |g_j(\lambda)|^2\:\Var{Y_j}\right)
\sup_{1\le j\le J}
\frac{\left|g_j(\lambda) - \sum_{n=1}^N a_n(\lambda) g_j(\lambda_n^N)\right|^2}{|g_j(\lambda)|^2} \,.
\end{align}

To get~\eqref{eq:rbcv},
we now explain how to choose the $N$ coefficients $a_n(\lambda)$, $1\le n\le N$,
for each $\lambda\in\Lambda$ when $\lambda_n^N\in\Lambda$, $n=1,\ldots,N$ is given, 
and then how to choose those $N$ parameter values $\lambda_n^N\in\Lambda$, $n=1,\ldots,N$.

Assume the $N$ parameter values $\lambda_n^N\in\Lambda$, $n=1,\ldots,N$, are given,
with $\lambda_0^N=\lambda_{\min}$, $\lambda_N^N=\lambda_{\max}$
and $\lambda_n^N\le\lambda_{n+1}^N$, $n=0,\ldots,N-1$.
Then, for a given $M\in\{2,\ldots,N\}$ (to be determined later on) 
and for all $\lambda\in\Lambda$,
it is possible to choose $1\le M_0(\lambda)\le N$ such that 
$\lambda^N_{M_0(\lambda)}\le\lambda\le\lambda^N_{M_0(\lambda)+M-1}$.
Only the $M$ coefficients corresponding to the $M$ contiguous parameter values above 
are taken non zero, such that $\forall\lambda\in\Lambda$~: 
$${a}_m(\lambda)\neq 0 \Leftrightarrow M_0(\lambda)\le m\le M_0(\lambda)+M-1 \,,$$
and are more specifically chosen as ${a}_m(\lambda)=P^\lambda_m(\tau_\Lambda(\lambda))$ where
$P^\lambda_m$ are polynomials of degree $M-1$,
such that, for all $M_0(\lambda)\le m,k\le M_0(\lambda)+M-1$,
$P^\lambda_m(\tau_\Lambda(\lambda_k))=\delta_{mk}$. 
The polynomial function $P^\lambda_m$ is the Lagrange interpolant defined on 
$[\tau_\Lambda(\lambda_{M_0(\lambda)}),\tau_\Lambda(\lambda_{M_0(\lambda)+M-1})]$, 
taking value $1$ at 
$\tau_\Lambda(\lambda_{m})$
and $0$ at $\tau_\Lambda(\lambda_{k})$, $k\neq m$.
We will also need a function $d(\lambda)= |\tau(\lambda_{M_0(\lambda)})-\tau(\lambda_{M_0(\lambda)+M-1})|$.
Using 
a Taylor-Lagrange formula for $g_j\circ\tau^{-1}$,
we have (for some $0\le \eta\le 1$)~:
\begin{equation*}
g_j(\lambda) - \sum_{n=1}^N a^\lambda_n(\lambda) g_j(\lambda_n)
=
\frac{d(\lambda)^{M}}{M!}
(g_j\circ\tau^{-1})^{(M)}
\left(
\eta\tau(\lambda_{M_0(\lambda)}^N)
+(1-\eta)\tau(\lambda_{M_0(\lambda)+M-1}^N)
\right)
\,.
\end{equation*}
Then, using~\eqref{eq:gassumption}
and the fact that $\Var{Z^\lambda} = \sum_{j=1}^J |g_j(\lambda)|^2\:\Var{Y_j}$,
there exists a constant $C>0$ (independent of $\Lambda$ and $J$) such that:
\begin{equation}\label{eq:bound2}
\Var{Z^\lambda-Y_N} \le 
\Var{Z^\lambda}\, \left(C\:d(\lambda)\right)^{2 M}
\,,\ \forall\lambda\in\Lambda\,.
\end{equation}

Finally, to get the result, we now choose a $\tau_\Lambda$-equidistributed parameter sample~:
$$ \tau_\Lambda(\lambda_n^N) = \tau_\Lambda(\lambda_{\min}) + \frac{n-1}{N-1} \left(\tau_\Lambda(\lambda_{\max})-\tau_\Lambda(\lambda_{\min})\right) 
\,,\ n=1,\ldots,N\,. $$
Then, $d(\lambda)=\frac{M-1}{N-1} \left(\tau_\Lambda(\lambda_{\max})-\tau_\Lambda(\lambda_{\min})\right)$ 
does not depend on $\lambda$.
Minimizing $(C\,d)^d$ as a function of $d\in(0,\frac1C)$, we choose $d(\lambda)=\frac1{e\:C}$,
and the choice 
$M=1+\lfloor \frac1{e\:C} \frac{N-1}{\tau_\Lambda(\lambda_{\max})-\tau_\Lambda(\lambda_{\min})} \rfloor$ (where $\lfloor x\rfloor$ denotes the integer part of a real number $x\in\R$) 
finishes the proof provided 
$N \ge N_0 \equiv 1+ \lfloor C\:e\:
\left(\tau_\Lambda(\lambda_{\max})-\tau_\Lambda(\lambda_{\min})\right) \rfloor$. $\ \square$

\end{document}